\documentclass[11pt]{amsart}

\usepackage{amsmath,amsthm,amsfonts,amssymb,bm,graphicx}

\usepackage%[backref=page]
{hyperref}
\hypersetup{colorlinks=true,citecolor=blue,linkcolor=blue,urlcolor=blue,
pdfstartview=FitH, pdfauthor=Valentin Blomer,pdftitle=On the 4-norm of an automorphic form}

\theoremstyle{plain}
\newtheorem{theorem}{Theorem}
\newtheorem{cor}{Corollary}
\newtheorem{lemma}{Lemma}

\numberwithin{equation}{section}

\theoremstyle{definition}

\renewcommand{\geq}{\geqslant}
\renewcommand{\leq}{\leqslant}

\makeatletter
\DeclareRobustCommand\widecheck[1]{{\mathpalette\@widecheck{#1}}}
\def\@widecheck#1#2{%
    \setbox\z@\hbox{\m@th$#1#2$}%
    \setbox\tw@\hbox{\m@th$#1%
       \widehat{%
          \vrule\@width\z@\@height\ht\z@
          \vrule\@height\z@\@width\wd\z@}$}%
    \dp\tw@-\ht\z@
    \@tempdima\ht\z@ \advance\@tempdima2\ht\tw@ \divide\@tempdima\thr@@
    \setbox\tw@\hbox{%
       \raise\@tempdima\hbox{\scalebox{1}[-1]{\lower\@tempdima\box
\tw@}}}%
    {\ooalign{\box\tw@ \cr \box\z@}}}
\makeatother

\begin{document}

\author{Valentin Blomer}
\address{Mathematisches Institut, Bunsenstr. 3-5, 37073 G\"ottingen, Germany} \email{blomer@uni-math.gwdg.de}

\title{On the 4-norm of an automorphic form}
 
\thanks{Author supported   by a Volkswagen Lichtenberg Fellowship and a Starting Grant of the European Research Council.   }

\keywords{$L^p$-norm, triple product $L$-functions, character sums}

\begin{abstract} We prove the optimal upper bound  $\sum_f \| f \|_4^4 \ll q^{\varepsilon}$ where $f$ runs over an orthonormal basis of Maa{\ss} cusp forms of prime level $q$ and bounded spectral parameter. 
\end{abstract}

\subjclass[2010]{11F12, 11F67 }

\setcounter{tocdepth}{2}  \maketitle %\tableofcontents

\section{Introduction}

Bounding $L^p$-norms of functions on a Riemannian surface (for $2 < p \leq \infty$) can be regarded as a weak type of equidistribution statement.  The situation is particularly interesting for manifolds with additional symmetries, such as a commutative algebra of Hecke operators commuting with the Laplacian, among other things because one can consider joint eigenfunctions which may rule out high multiplicity of eigenspaces. Often the underlying manifold is kept fixed, and one searches for bounds in terms of the Laplacian eigenvalue $\lambda$ as $\lambda \rightarrow \infty$. Here the first breakthrough for an arithmetic hyperbolic surface in the case $p=\infty$ has been obtained by Iwaniec and Sarnak \cite{IS}. 

In this article we change the point of view and keep the spectral data fixed, but study instead  the dependence on the manifold. We are interested in the 4-norm of a Maa{\ss} form on a hyperbolic surface $X_0(q) := \Gamma_0(q)\backslash \Bbb{H}$ where $q$ is a large prime. Equipped with the inner product
\begin{equation}\label{inner}
 \langle f, g\rangle = \int_{X_0(q)} f(z) \bar{g}(z) \frac{dx \, dy}{y^2},
\end{equation}
the space $X_0(q)$ has volume
\begin{equation}\label{vol}
  V(q) := \text{vol}(X_0(q))= \frac{\pi}{3} (q+1). 
\end{equation}

The 4-norm is a particularly interesting object because it is connected to triple product $L$-functions; 
by Watson's formula one has an equality roughly of the type
\begin{equation}\label{Watson}
  \| f \|_4^4 \approx  \frac{1}{q^2} \sum_{t_g \ll 1} L(1/2, f \times \bar{f} \times g)
\end{equation}
where the sum runs over an orthonormal basis of Hecke eigenforms $g$ of level $q$ with bounded spectral parameter $t_g$ (see \eqref{KimS} below).  By Weyl's law, the sum on the right hand side of \eqref{Watson} has $O(q)$ terms, so the Lindel\"of hypothesis for the $L$-functions on the right hand side of \eqref{Watson} would imply $\| f \|_4 \ll q^{-1/4 + \varepsilon}$, and  this is best possible by \eqref{vol}.  

The same  type of period formula is also the starting point for bounding the 4-norm in the \emph{eigenvalue} aspect, and in this case Sarnak and Watson have a announced a complete solution (possibly under the Ramanujan conjecture). Often in the theory of automorphic forms the archimedean and non-archimedean parameters behave, at least on a large scale, similarly. In the problem of bounding 4-norms, however, the spectral, weight and level aspect behave very differently: in  spectral  $t$ aspect, Watson's formula produces a sum of length $t^2$ of central $L$-values of conductor $t^8$, while in the weight $k$  aspect, Watson's formula produces a sum of length $k$ and conductor $k^6$ which is much harder to treat.  The level  aspect, that we focus on here, is also more difficult than the spectral aspect: %of Sarnak and Watson: 
the conductor of the   $L$-functions in \eqref{Watson} is $q^5$, so  again there is little hope to establish a Lindel\"of-type bound unconditionally with present technology. The aim of this article is to confirm this bound on average over Maa{\ss} forms $f$ of level $q$:
\begin{theorem}  Fix any real number $T > 1$ and any $\varepsilon > 0$. Then
\begin{equation}\label{mainthm}
  \sum_{t_f \leq T} \| f \|^4_4 \ll_{T, \varepsilon} q^{ \varepsilon}
\end{equation}
where the sum runs over an orthonormal basis of Maa{\ss} cusp forms of prime level $q$ and (fixed) spectral parameter $t_f \leq T$. 
\end{theorem}

Up to the power $q^{\varepsilon}$ Theorem 1 is best possible. For an individual form $f$, we have the trivial bound
\begin{displaymath}
  \| f \|_4 \leq \|f \|_2^{1/2} \| f \|_{\infty}^{1/2}. 
\end{displaymath}
Non-trivial bounds for $\| f \|_{\infty}$ have been obtained first in \cite{BH}, and the strongest result \cite{HT} implies
\begin{displaymath}
  \| f \|_{4} \leq q^{-1/12+\varepsilon}
\end{displaymath}
for an $L^2$-normalized Maa{\ss} form. It seems to be very hard to improve this on the basis of \eqref{Watson}. 
%Using \eqref{Watson} this could be improved, probably to $q^{-1/16+\varepsilon}$. 
Theorem 1 implies immediately the best possible bound $\| f \|_4 \ll q^{-1/4 +\varepsilon}$ for almost all $f$:
\begin{cor} For any $\delta > 0$ the bound $\| f \|_4 \ll q^{-1/4 + \delta}$ holds for all but $O(q^{1-4\delta+\varepsilon})$ of all Maa{\ss} forms $f$ occurring in the sum in \eqref{mainthm}. 
\end{cor}

The bound of Theorem 1 holds also for holomorphic cusp forms $f\in S_k(q)$ of any (fixed)  weight $k \geq 2$ and large prime level $q$. If $k$ is sufficiently large, one can use the Petersson formula instead of the Kuznetsov formula. For small $k$, one can  embed the holomorphic spectrum of weight $k$ into the Maa{\ss} spectrum of weight $k$ and use an appropriate weight $k$ Kuznetsov formula, see \cite{DFI}. \\

We remark on the side that the proof of Theorem 1 is  dependent on moderately strong bounds towards the Ramanujan conjecture. Any bound $|\Re \mu_{\pi}(q, i)| \leq 1/2 - \delta$ for the Langlands parameters associated to a cuspidal representation on $\pi$ on $GL_2$  at the (unramified) place $v = q$, as well as the archimedean bound $|\Re \mu_{\pi}(\infty, i)| \leq 1/2 - \delta$ for $\pi$ on $GL_2$ and  $GL_6$ suffices. Alternatively, if one prefers to stay entirely in $GL_2$,  then  $|\Re \mu_{\pi}(\infty, i)| \leq 1/6 - \delta$ for $\pi$ on $GL_2$ suffices.  In addition, we use several deep facts such as the automorphy of $GL_2 \times GL_3$ $L$-functions \cite{KSh}, non-negativity of central values, and of course Watson's   formula.  

%The estimation of the various off-diagonal terms is unchanged, the only difference is that Petersson's formula and Kuznetsov's formula have to be used \emph{together} as in \cite[(3.17) - (3.21)]{BHM}.   \\

It follows from the period formula \eqref{Watson} that the sum on the left hand side of Theorem 1  is roughly given by 
\begin{equation}\label{triple}
  q^{-2} \sum_{t_f, t_g \ll 1} L(1/2, f  \times \bar{f} \times g). 
\end{equation}
The (seemingly) similar average $\sum_{f, g} L(1/2, f \times g \times h)$ for $f, g, h \in S_2(q)$ holomorphic forms of weight 2 and level $q$  has been studied in \cite{FW}, also on the basis of triple product identities, but using entirely different techniques. 

There is another period formula in which the triple product $L$-functions in \eqref{triple} occur, %\footnote{I would like to than M.\ Young for pointing this out.}, 
namely as restrictions of certain Yoshida lifts. Given two holomorphic cuspidal Hecke forms $h_1, h_2 \in S_2(q)$ (more general assumptions are possible), one can  define  the (second) Yoshida lift $Y^{(2)}(h_1, h_2)$ which is a Siegel modular form of degree 2 and weight 2. When restricted to the diagonal $\left(\begin{smallmatrix} z_1 & \\ & z_2\end{smallmatrix}\right)$, it is a modular form of weight 2 both in $z_1$ and $z_2$, and hence
\begin{displaymath}
\begin{split}
 & Y^{(2)}(h_1, h_2)\left(\begin{smallmatrix} z_1 & \\ & z_2\end{smallmatrix}\right) =  \sum_{f_1, f_2 \in S_2(q)} c(f_1, f_2) f_1(z_1) f_2(z_2), \\
  & c(f_1, f_2) =  \int_{X_0(q)}\int_{X_0(q)}  Y^{(2)}(h_1, h_2)\left(\begin{smallmatrix} z_1 & \\ & z_2\end{smallmatrix}\right)f(z_1)f(z_2) \frac{dx_1\, dy_1}{y_1^2} \frac{dx_2 \, dy_2}{y_2^2}.
 \end{split} 
\end{displaymath}
A special case of a beautiful formula of B\"ocherer, Furusawa and Schulze-Pillot \cite[Corollary 2.7b]{BFSP} shows that for $h_1 = h_2=h$ and $f_1 = f_2 = f$ the coefficient $c(f, f)$ is proportional to the central $L$-value $L(1/2, f \times f \times h)$. The quantity estimated in Theorem 1 can then be interpreted as the trace of the matrix $(c(f_1, f_2))$, averaged over cusp forms $h$. \\

The paper is organized as follows: Sections 2 - 4 and 7 contain auxiliary material on automorphic forms, $L$-functions, character sums and integrals of Bessel functions. In particular we provide computations with oldforms,  newforms and Eisenstein series, a special type of approximate functional equation for the $L$-functions in question, and  bounds for certain complete exponential sums and oscillating integrals  that occur later in the analysis. Section 5 contains the main transformation from the average of 4-norms into smooth sums over products of Kloosterman sums that are estimated in Section 8. \\

I would like to thank M.\ Young and R.\ Schulze-Pillot for useful comments.

%We conclude the introduction by giving an informal description of the argument that underlies Theorem 1. An approximate functional equation

\section{Fourier expansions}

The spectrum of $L^2(X_0(q))$ consists of the constant function, Maa{\ss} forms, and Eisenstein series  $E_{\infty}(., 1/2 + it)$, $E_0(., 1/2 + it)$ for $t \in \Bbb{R}$, corresponding to the two ($\Gamma_0(q)$-equivalence classes of) cusps $\mathfrak{a} = \infty, 0$.  For any Maa{\ss} form $g$ we denote by %More precisely, we have
%\begin{displaymath}
%  E_{\infty}(z, s) = \frac{1}{2} \sum_{\Gamma_{\infty} \backslsh\Gamma_0(q)} \Im(\gamma z)^s = y^s + \end{displaymath}
\begin{equation}\label{KimS}
  t_g = \sqrt{ \lambda_g - 1/4} \in \mathcal{T} := \Bbb{R} \cup (-1/2, 1/2)i
\end{equation}  
its spectral parameter. % (cf.\ \cite{KS}). 

Let $\mathcal{B}_q$ be an orthonormal basis of cuspidal Hecke-Maa{\ss} newforms for $\Gamma_0(q)$.  Let $\mathcal{B}_1$ be a  basis of Hecke-Maa{\ss} cusp forms for $SL_2(\Bbb{Z})$ that is orthonormal with respect to the inner product \eqref{inner}.  In particular, for $g \in \mathcal{B}_1$ one has trivially 
\begin{equation}\label{upper}
  \| g \|_{\infty} \ll_{t_g} q^{-1/2}% t_g^{1/4}
\end{equation}  
by \eqref{vol}. The implied constant depends polynomially on $t_g$, for instance $(1+|t_g|)^{1/4}$ is admissible.  

For any such Hecke-Maa{\ss} cusp form $g$ in $\mathcal{B}_q$ or $\mathcal{B}_1$ we write $\lambda_g(n)$ for the $n$-th Hecke eigenvalue, and % and  $Q_{g}$ for its level, i.e. $Q_g = 1$ if $g \in \mathcal{B}_1$ and $Q_g = q$ if $g \in \mathcal{B}_q$. 
 we put $\delta_g = 0$ if $g$ is even and $\delta_g = 1$ if $g$ is odd. 

Newforms $g \in \mathcal{B}_q$ have two properties that we need later: they are eigenfunctions of the Fricke involution $z \mapsto -1/(qz)$, and one has 
\begin{equation}\label{small}
  \lambda_g(q) = \pm q^{-1/2}.
\end{equation}    
By Weyl's law we have   
\begin{equation}\label{weyl}
 \#\{ g \in \mathcal{B}_1 \mid t_{g} \leq T\}  \ll T^2, \quad \#\{g \in \mathcal{B}_q \mid t_g \leq T\} \ll qT^2.
\end{equation}

For $g \in \mathcal{B}_1$  define 
\begin{displaymath}
  g_q(z) := \Bigl(1 - \frac{q \lambda_g^2( q)}{(q+1)^2}\Bigr)^{-1/2} \Bigl(g(qz) - \frac{\lambda_g( q) q^{1/2}}{q+1} g(z)\Bigr). 
\end{displaymath}
By \cite[Proposition 2.6]{ILS}, $g$ and $g_q$ have the same norm and are orthogonal to each other.  We conclude that 
\begin{displaymath}
 \mathcal{B} :=   \mathcal{B}_q \cup \mathcal{B}_1  \cup \mathcal{B}_1^{\ast}, \quad \mathcal{B}_1^{\ast} := \{g_q \mid g \in \mathcal{B}_1\},
\end{displaymath}
is an orthonormal basis (with respect to \eqref{inner}) of the non-trivial cuspidal spectrum of $L^2(X_0(q))$.  

Let 
\begin{equation}\label{four}
  g(z) = \rho_g(1) \sum_{n \not= 0} \lambda_g(n) \sqrt{y} K_{it_g}(2 \pi n y)  e(nx)
\end{equation} 
be the Fourier expansion of some $g$ in $\mathcal{B}_1$ or $\mathcal{B}_q$ where $\lambda_g(-n) = \pm \lambda_g(n)$ depending on whether $g$ is even or odd.  We have the Rankin-Selberg bound
\begin{equation}\label{ranksel1}
  \sum_{n \leq x} |\lambda_g(n)|^2 \ll x (q (1+|t_g|) x)^{\varepsilon}.
\end{equation}
and the individual bound 
\begin{equation}\label{KimS1}
  \lambda_g(n) \ll n^{1/2 - \delta} 
\end{equation}
for some $\delta > 0$. Since $\underset{s=1}{\text{res}} E_{\infty}(z, s) = V(q)^{-1}$, we can compute 
\begin{equation}\label{unfolding}
\begin{split}
 1 & =  \| g \|_2^2  = V(q)  \underset{s=1}{\text{res}} \int_{X_0(q)} |g(z)|^2 E_{\infty}(z, s) \frac{dx\, dy}{y^2} \\
 &= |\rho_g(1)|^2 V(q) \underset{s=1}{\text{res}} \sum_{n \not=0} \frac{|\lambda_g(n)|^2}{n^s} \int_0^{\infty} K_{it_g}(2 \pi y)^2 y^{s} \frac{dy}{y} \\
 & = |\rho_g(1)|^2 V(q) 2 L(1, \text{Ad}^2 g) \frac{\pi}{8 \cosh(\pi t_g)} \begin{cases}\zeta^{(q)}(2)^{-1}, & g \in \mathcal{B}_q,\\ \zeta(2)^{-1}, & g \in \mathcal{B}_1.\end{cases}
\end{split}  
\end{equation} 
%where the superscript $^{(Q_g)}$ denotes that the Euler factor at $Q_g$ has been removed (this condition is empty if $Q_g = 1$).  
We conclude 
\begin{equation}\label{rho}
 | \rho_g(1)|  =  \left(\frac{2 \cosh(\pi t_g)}{  L(1, \text{Ad}^2g)}\right)^{1/2}\times  \begin{cases} (q+1)^{-1/2}, & g \in \mathcal{B}_1,\\
 (q')^{-1/2}, & g \in \mathcal{B}_q\end{cases}
\end{equation}
where
\begin{equation}\label{deftildeq}
  q' := \frac{q^2}{q-1} \asymp q.
\end{equation}
Let $g \in \mathcal{B}_1$ and let us define
\begin{displaymath}
  \lambda^{\ast}_g(n) := \Bigl(1 - \frac{q \lambda_g^2( q)}{(q+1)^2}\Bigr)^{-1/2} \Bigl(q^{1/2} \lambda_g\Bigl(\frac{n}{q}\Bigr) -  \frac{\lambda_g( q) q^{1/2}}{q+1} \lambda_g(n)\Bigr)
 \end{displaymath}
 with the convention $\lambda_g(x) = 0$ for $x \in \Bbb{Q} \setminus \Bbb{Z}$.  Then each $g \in \mathcal{B}_1^{\ast}$ has a Fourier expansion of type \eqref{four} with $\lambda_g^{\ast}(n)$ in place of  $\lambda_g(n)$. For $g \in \mathcal{B}_1$ and $q \nmid nm$ it follows that
  \begin{equation}\label{sum}
 \begin{split}
   \lambda_g(n)\bar{\lambda}_g(m) +& \lambda^{\ast}_g(n)\bar{\lambda}^{\ast}_g(m) = c_1(g, q)\lambda_g(n)\bar{\lambda}_g(m), \quad   c_1(g, q) =  \Bigl(1 - \frac{\lambda_g(q)^2 q}{(q+1)^2}\Bigr)^{-1} \asymp 1
   \end{split}
 \end{equation}
by \eqref{KimS1} and
\begin{equation}\label{sumq}
 \begin{split}
&  q^{1/2}\bigl( \lambda_g(qn)\bar{\lambda}_g(m) + \lambda^{\ast}_g(qn)\bar{\lambda}^{\ast}_g(m) \bigr)= c_2(g, q)\lambda_g(n)\bar{\lambda}_g(m), \\
   &  c_2(g, q) = q^{1/2} \lambda_g(q)\Bigl( 1 -\frac{q}{q+1}\Bigl(1 - \frac{\lambda^2_g(q)}{q+1} \Bigr)\Bigl( 1 -\frac{q\lambda_g(q)^2}{(q+1)^2}\Bigr)^{-1}\Bigr) \ll q^{-1/2} |\lambda_g(q)|  \ll 1. 
   \end{split}
 \end{equation}
The main point here is that even though for $g\in \mathcal{B}_1$ the formula \eqref{small} does not hold, an appropriate analogue is true if one combines the Fourier coefficients of $g$ and $g_q$. 
 
 Similar Fourier expansions hold for the Eisenstein series $E_{\mathfrak{a}}(z, s)$. Let
 \begin{displaymath}
   \eta(n, t) := \sum_{ad = |n|} (a/d)^{it}.
 \end{displaymath}
 Then
 \begin{displaymath}
  E_{\mathfrak{a}}(z, 1/2+it)   = \delta_{\mathfrak{a} = \infty} y^{1/2+it} + \phi_{\mathfrak{a}}(1/2 + it) y^{1/2-it} + \rho_{\mathfrak{a}}(1, t) \sum_{n \not= 0} \eta_{\mathfrak{a}}(n, t) \sqrt{y} K_{it}(2 \pi n y)  e(nx)
 \end{displaymath}
 where $\phi_{\mathfrak{a}}(s)$ is a meromorphic function that we do not need to specify, and  (see \cite[(3.25)]{CI})
 \begin{equation}\label{defeta}
 \begin{split}
 &  |\rho_{\mathfrak{a}}(1, t)| = \left(\frac{4 \cosh(\pi t)}{q|\zeta^{(q)}(1 + 2it)|}\right)^{1/2},\\
   & \eta_{\infty}(n, t) = \frac{\eta(n, t)}{q^{1/2+ it}} - q^{1/2} \eta(n/q, t), \quad \eta_0(n, t) = \eta(n, t) - q^{-it} \eta(n/q, t)
   \end{split}
 \end{equation}
 with the above convention that $\eta(x, t) = 0$ for $x \in \Bbb{Q} \setminus \Bbb{Z}$. For $q \nmid mn$ it follows that 
 \begin{equation}\label{eisen}
   \eta_{\infty}(n, t) \eta_{\infty}(m, -t) + \eta_0(n, t) \eta_0(m, -t) = \Bigl(1 + \frac{1}{q}\Bigr)\eta(n, t)\eta(m, -t)
 \end{equation}
and 
 \begin{equation}\label{eisen1}
  q^{1/2}\bigl( \eta_{\infty}(qn, t) \eta_{\infty}(m, -t) + \eta_0(qn, t) \eta_0(m, -t)\bigr) =  \frac{\eta(q, t)}{q^{1/2}} \eta(n, t) \eta(m, -t). 
 \end{equation}\\
 
 One of the most important tools is the Kuznetsov formula. Let $n, m \in \Bbb{Z}  $ be coprime to $q$ (in particular non-zero), and let $h$ be an even holomorphic function in $|\Re t| < 3/4$ such that $h(t) \ll (1+ |t|)^{-3}$. Then the Kuznetsov formula \cite[p.\ 409]{IK} together with the previous calculations  \eqref{rho}, \eqref{defeta}, \eqref{sum} and \eqref{eisen}  implies that 
\begin{equation}\label{kuz1}
\begin{split}
& 2 \sum_{g \in \mathcal{B}_1} \frac{c_1(g, q)\lambda_g(n) \bar{\lambda}_g(m) }{(q+1) L(1, \text{Ad}^2 g)} h(t_g) + 2 \sum_{g \in \mathcal{B}_q} \frac{\lambda_g(n) \bar{\lambda}_g(m) }{q' L(1, \text{Ad}^2 g)} h(t_g) +   \int_{\Bbb{R}} \frac{\eta(n, t)\eta(m, -t)}{q'' |\zeta^{(q)}(1+2 i t)|^2} h(t) \frac{dt}{\pi} \\
&= \delta_{n, m} \int_{0}^{\infty} h(t) \frac{ d^{\ast}t}{\pi^2}  + \sum_{q \mid c} \frac{1}{c}S(n, m, c)  \int _0^{\infty} \mathcal{J}^{\pm}\Bigl(\frac{\sqrt{|nm|}}{c}, t\Bigr) h(t)  \frac{d^{\ast} t}{\pi} %h^{\pm}\left(\frac{\sqrt{|mn|}}{c}\right)
\end{split}
\end{equation}
%\begin{displaymath}
%\begin{split}
%& 2 \sum_{g \in \mathcal{B}_1} \frac{\lambda_g(n) \lambda_g(m) + \lambda_g^{\ast}(n)\lambda_g^{\ast}(m)}{(q+1) L(1, \text{Ad}^2 g)} h(t_g) + 2 \sum_{g \in \mathcal{B}_q} \frac{\lambda_g(n) \lambda_g(m) }{\tilde{q} L(1, \text{Ad}^2 g)} h(t_g) +  \sum_{\mathfrak{a}} \int_{\Bbb{R}} \frac{\eta_{\mathfrak{a}}(n; t)\eta_{\mathfrak{a}}(m; t)}{q|\zeta^{(q)}(1+2 i t)|^2} h(t) \frac{dt}{\pi} \\
%&= \delta_{n, m} \int_{0}^{\infty} h(t) \frac{ d^{\ast}t}{\pi^2}  + \sum_{q \mid c} \frac{1}{c}S(n, m, c)  \int _0^{\infty} \mathcal{J}^{\pm}(x, t) h(t)  \frac{d^{\ast} t}{\pi} %h^{\pm}\left(\frac{\sqrt{|mn|}}{c}\right)
%\end{split}
%\end{displaymath}
where $\pm = \text{sgn}(mn)$, $d^{\ast} t = t \tanh(\pi t) dt$, $q' = q^2/(q-1)$ as in \eqref{deftildeq}, $q'' = q^2/(q+1)$  and
%\begin{displaymath}
 % h^{\pm}(x) =  \int _0^{\infty} \mathcal{J}^{\pm}(x, t) h(t)  \frac{d^{\ast} t}{\pi}. \end{displaymath}
 % with
  \begin{displaymath}
    \mathcal{J}^{\pm}(x, t) =   \begin{cases} \displaystyle  \frac{2  i}{ \sinh(\pi t)} (J_{2 it}(4 \pi x) - J_{-2it}(4 \pi x)),\\ \displaystyle \frac{2  i}{ \sinh(\pi t)} (I_{2 it}(4 \pi x) - I_{-2it}(4 \pi x)) = \frac{4}{\pi} K_{2 i t}(4 \pi x) \cosh(\pi t). \end{cases}
  \end{displaymath}
Similarly, for $q\nmid nm$ we obtain by \eqref{sumq} and \eqref{eisen1} instead of \eqref{sum}  and \eqref{eisen} that 
\begin{equation}\label{kuz2}
\begin{split}
& 2 \sum_{g \in \mathcal{B}_1} \frac{c_2(g, q)\lambda_g(n) \bar{\lambda}_g(m) }{(q+1) L(1, \text{Ad}^2 g)} h(t_g) + 2 q^{1/2} \sum_{g \in \mathcal{B}_q} \frac{\lambda_g(qn) \bar{\lambda}_g(m) }{q' L(1, \text{Ad}^2 g)} h(t_g) +   \int_{\Bbb{R}} \frac{\eta(qn, t)\eta(m, -t)}{q^{3/2} |\zeta^{(q)}(1+2 i t)|^2} h(t) \frac{dt}{\pi} \\
&=    q^{1/2} \sum_{q \mid c} \frac{1}{c}S(qn, m, c)  \int _0^{\infty} \mathcal{J}^{\pm}\Bigl(\frac{\sqrt{|qnm|}}{c}, t\Bigr) h(t)  \frac{d^{\ast} t}{\pi} . %h^{\pm}\left(\frac{\sqrt{|mn|}}{c}\right)
\end{split}
\end{equation}

Let
\begin{displaymath}
  \mathcal{J}^{0}(x, t) := \frac{1}{2} \bigl(   \mathcal{J}^{+}(x, t) +  \mathcal{J}^{-}(x, t)\bigr).
\end{displaymath}
Adding the   Kuznetsov formula for $nm > 0$ and $nm < 0$, we can single out even Maa{\ss} forms:  
\begin{equation}\label{kuz3}
\begin{split}
& 2 \sum_{\substack{g \in \mathcal{B}_1\\ g \text{ even}}} \frac{c_1(g, q) \lambda_g(n) \bar{\lambda}_g(m) }{(q+1) L(1, \text{Ad}^2 g)} h(t_g) + 2 \sum_{\substack{g \in \mathcal{B}_q\\ g \text{ even}}} \frac{\lambda_g(n) \bar{\lambda}_g(m) }{q' L(1, \text{Ad}^2 g)} h(t_g) +  \sum_{\mathfrak{a}} \int_{\Bbb{R}} \frac{\eta(n; t)\eta(m; -t)}{q''|\zeta^{(q)}(1+2 i t)|^2} h(t) \frac{dt}{\pi} \\
&= \delta_{n, m} \int_{0}^{\infty} h(t) \frac{ d^{\ast}t}{2\pi^2}  + \sum_{q \mid c} \frac{1}{c}S(n, m, c) \int_0^{\infty} \mathcal{J}^0\Bigl(\frac{\sqrt{nm}}{c}, t\Bigr) h(t) \frac{d^{\ast} t}{\pi}% \left(h^{+}\left(\frac{\sqrt{mn}}{c}\right)+ h^{-}\left(\frac{\sqrt{mn}}{c}\right)\right)
\end{split}
\end{equation}  
for $m, n \in \Bbb{N}$, $q\nmid nm$, as well as 
\begin{equation}\label{kuz4}
\begin{split}
& 2 \sum_{\substack{g \in \mathcal{B}_1\\ g \text{ even}}} \frac{c_2(g, q)\lambda_g(n) \bar{\lambda}_g(m) }{(q+1) L(1, \text{Ad}^2 g)} h(t_g) + 2 q^{1/2} \sum_{\substack{g \in \mathcal{B}_q\\ g \text{ even}}} \frac{\lambda_g(qn) \bar{\lambda}_g(m) }{q' L(1, \text{Ad}^2 g)} h(t_g) +   \int_{\Bbb{R}} \frac{\eta(qn, t)\eta(m, -t)}{q^{3/2} |\zeta^{(q)}(1+2 i t)|^2} h(t) \frac{dt}{\pi} \\
&=    q^{1/2} \sum_{q \mid c} \frac{1}{c}S(qn, m, c)  \int _0^{\infty} \mathcal{J}^{\pm}\Bigl(\frac{\sqrt{qnm}}{c}, t\Bigr) h(t)  \frac{d^{\ast} t}{\pi} . %h^{\pm}\left(\frac{\sqrt{|mn|}}{c}\right)
\end{split}
\end{equation}
We will need all 4 versions \eqref{kuz1} - \eqref{kuz4} in Section \ref{main}.

\section{Triple product $L$-functions} 

Let $f, g \in   \mathcal{B}_q$.  Then we can define the triple product $L$-function
\begin{displaymath}
  L(s, f\times \bar{f} \times g) = L(s, \text{Ad}^2 f \times g) L(s, g). 
\end{displaymath}
The local factors, root number and conductor have been computed in \cite[Section 4.1]{Wa}. 
   Let $\Gamma_{\Bbb{R}}(s) := \Gamma(s/2) \pi^{-s/2}$. Then
\begin{displaymath}
\begin{split}
& L_{\infty}(s, g) :=  \prod_{\pm} \Gamma_{\Bbb{R}}(s  \pm i t_g), \\
 &  \Lambda(s,   g) = L(s,   g) L_{\infty}(s, g)   = (-1)^{\delta_g} (-\lambda_g(q)q^{1/2}) q^{1/2 - s}\Lambda(1-s,   g).
\end{split} 
 \end{displaymath}
%where 
%\begin{displaymath}
%\epsilon_g = \begin{cases} (-1)^{\delta_g}, & g \in \mathcal{B}_1,\\
%(-1)^{\delta_g} (-\lambda_g(q)q^{1/2}), & g \in \mathcal{B}_q\end{cases}
%\end{displaymath}
%is the root number of $g$. 
Similarly, 
\begin{displaymath}
\begin{split}
&  L_{\infty}(s, \text{Ad}^2 f \times g) = \prod_{\pm} \prod_{\nu=-1}^1\Gamma_{\Bbb{R}}(s + 2 i \nu t_f \pm i t_g), \\
 & \Lambda(s, \text{Ad}^2 f \times g) = L(s, \text{Ad}^2f \times g)   L_{\infty}(s, \text{Ad}^2 f \times g)  = (-1)^{\delta_g}  (q^4)^{1/2 - s}\Lambda(1-s, \text{Ad}^2f \times g).
 \end{split}
 \end{displaymath}
By \cite[Theorem 1]{LR}, the central value $L(1/2, g)$ is non-negative. Moreover,  by \cite{JS} the self-dual lift $\text{Ad}^2 f \, (= \text{sym}^2 f)$ is orthogonal (its symmetric square $L$-function has a pole at $s=1$), hence  by  \cite[Theorem 1.1]{La}, the central value $\Lambda(1/2, \text{Ad}^2 f \times g)$ is non-negative, too, and an inspection of the archimedean $L$-factors shows that the same holds for  $L(1/2, \text{Ad}^2 f \times g)$, hence also $L(1/2, f\times \bar{f} \times g)$.   We note that $L(s, \text{Ad}^2 f \times g)$ and hence $L(s, f \times \bar{f} \times g)$ vanishes at $s=1/2$ if $g$ is odd. \\
%where
%\begin{displaymath}
 % Q_{f, g} = \begin{cases}1, & f, g \in \mathcal{B}_1, \\ q, & f \in \mathcal{B}_1, g \in \mathcal{B}_q,\\ q^2,& f \in \mathcal{B}_q, g \in \mathcal{B}_1,\\ q^4, & f, g \in \mathcal{B}_q. \end{cases}
%\end{displaymath}

The adjoint square lift of $f$ is a self-dual automorphic form on $GL_3$ with Fourier coefficients $A(m, k)$ satisfying $A(m, 1)  = \sum_{ab^2 = m} \lambda_f(a^2)$ for $q \nmid m$. Using Hecke relations, we can express all $A(m, k)$ in terms of the Hecke eigenvalues of $\lambda_f$ as follows: by M\"obius inversion and \cite[Theorem 6.4.11]{Go} we have
\begin{displaymath}
  A(m, k) = \sum_{d \mid (m, k)} \mu(d) A\left(\frac{m}{d}, 1\right) A\left(1, \frac{k}{d}\right)%= \sum_{d \mid (n, m)} \mu(d) \lambda_f\left(\left(\frac{n}{d}\right)^2\right) \lambda_f\left(\left(\frac{m}{d}\right)^2\right). 
\end{displaymath}
whenever $q \nmid mk$. Hence
\begin{equation}\label{coeff-q}
  L^{(q)}(s, \text{Ad}^2 f \times g) =  \sum_{q \nmid mk} \frac{A(m, k) \lambda_g(m)}{m^sk^{2s}} = \sum_{q \nmid dabmk}  \frac{\mu(d)\lambda_f(m^2)\lambda_f(k^2)\lambda_g(dma^2)}{m^sa^{2s}k^{2s}b^{4s}d^{3s}}.
\end{equation}
Using the explicit shape of the Euler factor at $q$ (see \cite{Wa}), we find
%For any $L$-function $L(s, \cdot)$ let $L^{(q)}(s, \cdot)$ denote the $L$-functions with the Euler factor at $q$ removed. Then
\begin{equation}\label{coeff}
\begin{split}
 &   L(s, \text{Ad}^2 f \times g) =  \left(1 - \frac{\lambda_g(q)}{q^s}\right)^{-1} \left(1 - \frac{\lambda_g(q)}{q^{s+1}}\right)^{-1}  L^{(q)}(s, \text{Ad}^2 f \times g)
 % \sum_{q \nmid n m} \frac{A(n, m) \lambda_g(n)}{n^sm^{2s}}\\
%    & = \sum_{k=0}^{\infty} \frac{\lambda_g(q)^k (1-q^{-k-1})}{q^{ks}(1-q^{-1})} \sum_{q \nmid dabnm}  \frac{\mu(d)\lambda_f(n^2)\lambda_f(m^2)\lambda_g(dna^2)}{n^sa^{2s}m^{2s}b^{4s}d^{3s}}
 =: \sum_{m} \frac{\lambda_{\text{Ad}^2 f \times g}(m)}{m^{s}},
 \end{split}    
\end{equation}
say. Note that by \eqref{small} the coefficients divisible by $q$ are small. Since these are purely formal computations with local Euler factors, \eqref{coeff-q} holds also  for $f$ and/or $g$ in $\mathcal{B}_1$, and analogous formulas hold for Eisenstein series:
\begin{displaymath}
\begin{split}
 &| L^{(q)}(s + it, \text{Ad}^2 f)|^2  = \sum_{q \nmid dabnm}  \frac{\mu(d)\lambda_f(m^2)\lambda_f(k^2)\eta(dma^2, t)}{m^sa^{2s}k^{2s}b^{4s}d^{3s}}, \\
 & | L^{(q)}(s + it, g)|^2 L^{(q)}(s, g)  = \sum_{q \nmid dabmk}  \frac{\mu(d)\eta(m^2, t)\eta(k^2, t)\lambda_g(dma^2)}{m^sa^{2s}k^{2s}b^{4s}d^{3s}}.
\end{split} 
\end{displaymath}
We have already seen that $L(s, f \times \bar{f} \times g)$ has conductor $q^5$ for $f, g \in \mathcal{B}_q$. If one of the factors has level one or is an Eisenstein series, the conductor drops; more precisely, all the $L$-functions
\begin{displaymath}
\begin{split}
 & L(s, f \times \bar{f} \times g), \quad |L(s + it, f \times \bar{f})|^2, \quad f \in \mathcal{B}_q, g \in \mathcal{B}_1,\\
  & L(s, f \times \bar{f} \times g), \quad |L(s + it, g)|^2 L(s, g)^2, \quad f \in \mathcal{B}_q, g \in \mathcal{B}_q
  \end{split}
\end{displaymath}
have conductor $q^4$. We will use this observation in Sections \ref{main} and  \ref{old}.

It is a deep result \cite{KSh} that $\text{Ad}^2 f \times g$ corresponds to an automorphic form on $GL(6)$. Hence the Rankin-Selberg $L$-function $L(s, (\text{Ad}^2 f \times g) \times (\text{Ad}^2 f \times g))$ satisfies the properties of \cite[Theorem 2]{Li}, and we have the upper bound
\begin{equation}\label{ranksel2}
  \sum_{m \leq x} |\lambda_{\text{Ad}^2 f \times g}(m)|^2 \ll x (q(1+|t_g|+|t_f|)x)^{\varepsilon}.
\end{equation}

We need a somewhat sophisticated and carefully designed approximate functional equation and borrow some ideas from \cite{Bl2}.  Let $A_1, A_2 \geq 10$ be integers and define
\begin{displaymath}
  G(u) = \left(\cos\frac{\pi u}{4A_1}\right)^{-100 A_1}, \quad 
%\end{displaymath}
%\begin{displaymath}
  G_1(u, t) := \prod_{\pm}  \prod_{\pm} \prod_{\ell=0}^{A_2} \left(\frac{1/2 \pm u  \pm it}{2}+\ell\right)  %\left(\frac{1/2 -u   \pm it }{2}+\ell\right),
\end{displaymath}
 and 
\begin{displaymath}
\begin{split}
  G_2(u, t_1, t_2) := \prod_{\pm} \prod_{\pm} \prod_{\nu=-1}^1\prod_{\ell=0}^{3A_2} &\left(\frac{1/2 \pm u \pm it_1 + 2 i \nu t_2}{2}+\ell\right). % \left(\frac{1/2 -u  \pm it_1 +2i\nu t_2}{2}+\ell\right).
%  \times   &  \left(\frac{1/2 +u + \delta \pm  it - \overline{\alpha_j}}{2}+\ell \right)  \left(\frac{1/2 -u + \delta \pm it - \overline{\alpha_j}}{2}+\ell \right).
 \end{split}
  \end{displaymath}
Clearly $G_1$ and $G_2$ are holomorphic and even in all variables, and $G$ is even and holomorphic in $|\Re u| < 2A_1$.  Moreover, for $t, t_1, t_2 \in \mathcal{T}$ we have
\begin{equation}\label{lowerG}
  G_1(0, t), \, G_2(0, t_1, t_2) \gg 1. % \asymp (1 + |t|)^{4(A_2+1)}, \quad G_2(0, t_1, t_2) \asymp \prod_{\nu=-1}^1 (1+|t_1 + 2 \nu t_2|)^{4(3A_2+1)}
\end{equation}
For this lower bound we either need that $|\Im  t_1|, |\Im t_2| \leq 1/6 - \delta$ or any nontrivial bound towards the Ramanujan conjecture for the infinite place of the $GL_6$ automorphic form $\text{Ad}^2 f \times g$. Both results are known \cite{KS, LRS}.  Let
\begin{displaymath}
  V_1(y; t) = \frac{1}{2\pi i} \int_{(2)} G(u) G_1(u, t) \prod_{\pm} \frac{\Gamma_{\Bbb{R}}(1/2 + u  \pm i t)}{\Gamma_{\Bbb{R}}(1/2  \pm i t)}  y^{-u} \frac{du}{u}
\end{displaymath}
and
\begin{displaymath}
 V_2(y; t_1, t_2) = \frac{1}{2\pi i} \int_{(2)} G(u) G_2(u, t_1, t_2) \prod_{\pm}\prod_{\nu=-1}^1 \frac{\Gamma_{\Bbb{R}}(1/2 + u  \pm i t_1 \pm 2 i \nu t_2)}{\Gamma_{\Bbb{R}}(1/2  \pm i t_1 \pm 2 i \nu t_2)}  y^{-u} \frac{du}{u}.
\end{displaymath}
 The weight functions $V_1, V_2$ have the following properties:
 \begin{lemma}\label{approx} a) The function $V_1(y; t)$ is smooth for $y > 0$ and holomorphic in $|\Im t| \leq 2 A_2$ and satisfies the uniform bound
 \begin{displaymath}
   y^jV_1^{(j)}(y; t) \ll (1+|t|)^{4(A_2+1)} \left(1 + \frac{y}{1+|t|}\right)^{-A_1} \ll (1+|t|)^{A_1 + 4(A_2+1)} (1+y)^{-A_1}
 \end{displaymath} 
 in this region for fixed $j \in \Bbb{N}_0$. Its Mellin transform with respect to the first variable, $\widehat{V}_1(u; t)$, is holomorphic in $\Re u \geq \varepsilon$ whenever $t \in \mathcal{T}$. %, except for  a simple pole at $u=0$ with residue $G_1(0, t)$.  
 In this region it satisfies the uniform bound
 \begin{displaymath}
\widehat{V}_1(u; t)  \ll_{\Re u, \varepsilon}  e^{-\Im |u|} (1+|t|)^{4(A_2+1) + \Re u}. %, \quad  |u| \geq 1/10. 
 \end{displaymath}
Moreover,
\begin{equation}\label{zero}
  \widehat{V}_1(1/2 \pm it, t) = 0.
\end{equation}

b) The function $V_2(y, t_1, t_2)$ is smooth in $y > 0$ and holomorphic in $|\Im t_1|, |\Im t_2|  \leq 2 A_2$ and satisfies the uniform bound
 \begin{displaymath}
 \begin{split}
  y^j V^{(j)}_2(y; t) &\ll (1+|t_1| +|t_2|)^{12(3A_2+1)} \left(1 + \frac{y}{(1+|t_1|+|t|_2)^3}\right)^{-A_1}\\
  & \ll (1+|t_1|+ |t_2|)^{3A_1 + 12(3A_2+1)} (1+y)^{-A_1}
  \end{split}
 \end{displaymath} 
in this region for fixed $j \in \Bbb{N}_0$. Its Mellin transform with respect to the first variable, $\widehat{V}_2(u; t_1, t_2)$, is holomorphic in $\Re u \geq \varepsilon$ whenever $t_1, t_2 \in \mathcal{T}$. %, except for  a simple pole at $u=0$ with residue $G_2(0, t_1, t_2)$. 
 In this region it satisfies the uniform bound
 \begin{displaymath}
\widehat{V}_2(u; t_1, t_2)  \ll_{\Re u, \varepsilon}  e^{-\Im |u|} (1+|t_1|+|t_2|)^{12(3A_2+1) + 3\Re u}. %, \quad |u| \geq 1/10.  
 \end{displaymath}
\end{lemma}

\textbf{Proof.} This follows easily from the definition of $G$, $G_1$, $G_2$ with appropriate contour shifts. \\

Let $g \in \mathcal{B}_q$ be even.   Then the usual technique (e.g.\ \cite[p.\ 98]{IK} or \cite[Section 2]{Bl2}) shows
  \begin{equation}\label{appr1}
   G_1(0, t_g) L(1/2, g) = (1-\lambda_g(q)q^{1/2}) \sum_{n} \frac{\lambda_g(n)}{n^{1/2}} V_1\left(\frac{n}{q^{1/2}}; t_g\right)% - \lambda_g(q)q^{1/2} \sum_{n} \frac{\lambda_g(n)}{n^{1/2}} V_1\left(\frac{n}{q^{1/2}}; t_g\right)
  \end{equation}
  and 
   \begin{equation}\label{appr2}
   G_2(0, t_g, t_f) L(1/2, \text{Ad}^2 f \times g ) = 2 \sum_{m} \frac{\lambda_{\text{Ad}^2f\times g}(m)}{m^{1/2}} V_2\left(\frac{m}{q^{2}}; t_g, t_f\right).  
  \end{equation}

\section{Character sums}

For future reference we state some useful results for character sums. We quote from \cite[Section 3]{Bl1}. For a positive or negative discriminant $D$ of a quadratic number field let $\chi_D = \left(\frac{D}{.}\right)$ be the associated Dirichlet character. Define 
\begin{displaymath}
  \epsilon_c := \left\{
    \begin{array}{ll}
       1, & c > 0,\\
        i, & c < 0,
     \end{array}\right.   
\end{displaymath}
and if $c=c_1c_2^2$ is odd and positive with $\mu^2(c_1)=1$ let 
\begin{displaymath}
  c^{\ast} := \chi_{-4}(c_1) c_1.
\end{displaymath}  
We need to evaluate the sum 
\begin{displaymath}
  G(d, h; q):= \sum_{x \, (q)} e\left(\frac{dx^2 + hx}{q}\right) 
\end{displaymath}
for integers $d \in \Bbb{Z} \setminus \{0\}$, $h \in \Bbb{Z}$, $q \in \Bbb{N}$. Clearly
\begin{displaymath}
  G(d, h; q) = \delta_{(d, q) \mid h} G(d/(d, q), h/(d, q); q/(d, q)),
\end{displaymath}
so it suffices to compute the sum for $(d, q) = 1$. We write $q = s2^{\alpha}$ with $s$ odd. Then we have \cite[Lemma 2]{Bl1}
\begin{equation}\label{gauss1}
G(d, h; q) = \begin{cases} 
  \sqrt{q} \epsilon_{q^{\ast}} \chi_{q^{\ast}}(d)e\left(\frac{-\overline{4d} h^2}{q}\right) , & \alpha = 0,\\
  \sqrt{2q} \epsilon_{s^{\ast}} \chi_{s^{\ast}}(2d) e\left(\frac{-\overline{8d}h^2}{s}\right), & \alpha = 1, h \text{ odd},\\
   \sqrt{q} \epsilon_{s^{\ast}} \chi_{s^{\ast}}(d) e\left(\frac{-\bar{d}(h')^2}{q}\right)(1+i \chi_{-4}(sd)), & \alpha \geq 2 \text{ even}, h  = 2h' \text{ even},\\
    \sqrt{q} \epsilon_{s^{\ast}} \chi_{8s^{\ast}}(d) e\left(\frac{-\bar{d}(h')^2}{q}\right)(1 +i\chi_{-4}(sd) ), &\alpha \geq 3 \text{ odd}, h = 2h' \text{ even},\\
    0, & \text{otherwise,}
\end{cases}
\end{equation}
whenever $(d, q) = 1$. 
If $\psi$ is a real character of conductor $s$, $q = ss_1s_2$ with $s_1 \mid s^{\infty}$ and $(s, s_2)= 1$, and $\Delta \in \Bbb{Z}$, then we have (\cite[(3.2)]{Bl1})
\begin{equation}\label{gauss2}
\begin{split}
 & \Bigl|\left. \sum_{d \, (q)}\right.^{\ast} \psi(d) e\left(\frac{d\Delta}{q}\right)\Bigr|   = \Bigl|\delta_{s_1 \mid \Delta} \psi\left(s_2\frac{\Delta}{s_1}\right)s_1 r_{s_2}(\Delta) \sqrt{s} \epsilon_{s^{\ast}}\Bigr| \leq \sqrt{qs_1s_2} 
\end{split}
\end{equation}
where $r_{q}(\Delta)$ is the Ramanujan sum.  

For a Schwartz class function $W$ we denote by $\check{W}$ its Fourier transform.

\begin{lemma}\label{char1}  Let $\gamma \in \Bbb{N}$, $\alpha \in \Bbb{Z}$, and let $W$ be a Schwartz class  function. Then 
\begin{displaymath}
  \sum_{n \in \Bbb{Z}} S(n, \alpha, \gamma) W(n) = \sum_{\substack{h_1 \in \Bbb{Z}\\ (h_1, \gamma) = 1}} e\left(- \frac{\alpha \bar{h}_1}{\gamma}\right) \check{W}\left(\frac{h_1}{\gamma}\right). 
\end{displaymath}
\end{lemma}

\textbf{Proof.} This is a direct consequence of the Poisson summation formula.\\

\begin{lemma}\label{char2} Let $c, r \in \Bbb{N}$, $\beta, \kappa, h  \in \Bbb{Z}$. % and assume $q \nmid cr$. 
Write $r = r_1 r_2^2$ with $\mu(r_1)^2 = 1$ and write $f = (c, r)$, $c = f \tilde{c}$, $r = f \tilde{r}$ with $(\tilde{r}, \tilde{c}) = 1$.  Then
\begin{displaymath}
    \Bigl|\sum_{m \, (q\tilde{c}\tilde{r}f)}  e\left(\frac{m \beta}{qc}\right) S(m^2, \kappa, qr)   e\left(\frac{mh}{q\tilde{c}\tilde{r}f}\right)\Bigr| \begin{cases} \leq  2qrr_2\tilde{c} , & \text{if }  \beta  \tilde{r} \equiv - h \, (\tilde{c}),\\
  = 0, & \text{otherwise.}\end{cases}
\end{displaymath}
\end{lemma}

\textbf{Proof.} Opening the Kloosterman sum, the exponential sum in question equals
\begin{displaymath}
  \underset{x \, (qr)}{\left.\sum \right.^{\ast}} e\left(\frac{\kappa \bar{x}}{qr}\right) G(x\tilde{c}, h+\beta \tilde{r}, q\tilde{c}\tilde{r}f) = \tilde{c} \delta_{\tilde{c} \mid h+\beta \tilde{r}}   \underset{x \, (qr)}{\left.\sum \right.^{\ast}} e\left(\frac{\kappa \bar{x}}{qr}\right) G(x, (h+\beta \tilde{r})/\tilde{c}, qr). 
\end{displaymath}
For notational simplicity let us write $\gamma := (h+\beta \tilde{r})/\tilde{c}$. We evaluate the Gau{\ss} sum using \eqref{gauss1}. To this end, we write $qr = s 2^a$ with $s$ odd and also recall $r = r_1r_2^2$ with $\mu(r_1)^2 = 1$. We distinguish several very similar cases. If $a = 0$, we obtain
\begin{displaymath}
 \sqrt{qr} \epsilon_{(qr)^{\ast}}  \tilde{c} \delta_{\tilde{c} \mid h+\beta \tilde{r}}   \underset{x \, (qr)}{\left.\sum \right.^{\ast}} e\left(\frac{\kappa \bar{x}}{qr}\right) \chi_{(qr)^{\ast}}(x) e\left(\frac{-\overline{4x}\gamma^2}{qr}\right)  \end{displaymath}
and the desired bound (without the factor 2) follows directly from \eqref{gauss2}. If $a = 1$, we obtain
\begin{displaymath}
\sqrt{2qr} \epsilon_{s^{\ast}}  \tilde{c} \delta_{\tilde{c} \mid h+\beta \tilde{r}}  \delta_{2 \nmid \gamma} \underset{x \, (qr)}{\left.\sum \right.^{\ast}} e\left(\frac{\kappa \bar{x}}{qr}\right) \chi_{s^{\ast}}(2x) e\left(\frac{-\overline{8x}\gamma^2}{s}\right).
   \end{displaymath}
   The $x$-sum equals in absolute value
   \begin{displaymath}
    \Bigl| \left.\sum_{x \, (s)}\right.^{\ast} e\left(\frac{\kappa \overline{2x}}{s}\right) \chi_{s^{\ast}}(2x) e\left(\frac{-\overline{8x}\gamma^2}{s}\right)\Bigr| \leq  \sqrt{s} r_2 = r_2 \sqrt{qr/2}
\end{displaymath}  
and the lemma follows again (without the factor 2). If $a \geq 2$ is even, we have
\begin{displaymath}
  \sqrt{qr} \epsilon_{s^{\ast}} \tilde{c} \delta_{\tilde{c} \mid h+\beta \tilde{r}}  \delta_{2 \mid \gamma} \underset{x \, (qr)}{\left.\sum \right.^{\ast}} e\left(\frac{\kappa \bar{x}}{qr}\right)\chi_{s^{\ast}}(x)  e\left(\frac{-\bar{x}(\gamma/2)^2}{qr}\right) (1 + i \chi_{-4}(sx))
\end{displaymath}
and the lemma follows from \eqref{gauss2}. The case $a \geq 3$ odd is identical. 

%\begin{lemma} Let $b, k, h_1, h_2 \in \Bbb{Z}$, $c, r \in \Bbb{N}$ and assume $q \nmid cr$. Write $r = r_1 r_2^2$ with $\mu(r_1)^2 = 1$ and write $d = (c, r)$, $c = d c'$, $r = d r'$ with $(r', c') = 1$ Then
%\begin{displaymath}
 % \sum_{n \, (qc)}  \sum_{m \, (qc'r'd)} S(n, mb, qc) S(m^2, k, qr) e\left(\frac{nh_1}{qc}\right) e\left(\frac{mh_2}{qc'r'd}\right) \begin{cases} \ll  q^2 r r_2 cc' , & \text{if } (h_1, qc) = 1 \text { and } \\
 % & b r' \equiv - h_1h_2 \, (c'),\\
 % = 0, & \text{otherwise.}\end{cases}
%\end{displaymath}
%\end{lemma}

%\textbf{Proof.} We open the two Kloosterman sums getting
%\begin{displaymath}
 %  \sum_{n \, (qc)}  \sum_{m \, (qcr/d)}\underset{x \, (qc)}{\left.\sum \right.^{\ast}}  \underset{y \, (qc'r'd)}{\left.\sum \right.^{\ast}} e\left(\frac{n x + mb \bar{x}}{qc}\right) e\left(\frac{m^2 y + k\bar{y}}{qr}\right)  e\left(\frac{nh_1}{qc}\right) e\left(\frac{mh_2}{qc'r'd}\right)
%\end{displaymath}
%We sum over $n$

\begin{lemma}\label{char3} 
 Let  $c, r \in \Bbb{N}$ and %assume $q \nmid cr$. 
 let  $\beta, \kappa  \in \Bbb{Z}$. Write $r = r_1 r_2^2$ with $\mu(r_1)^2 = 1$ and write $f = (c, r)$, $c = f \tilde{c}$, $r = f \tilde{r}$ with $(\tilde{r}, \tilde{c}) = 1$. Let $W$ be a Schwartz class function. Then
\begin{displaymath}
    \sum_{m \in \Bbb{Z}}  e\left(\frac{m \beta}{qc}\right) S(m^2, \kappa, qr)   W(m)  \ll  r_2 \sum_{\substack{h_2 \in \Bbb{Z}\\ h_2 \equiv -\beta \tilde{r} \, (\tilde{c})}} \Bigl|\check{W}\Bigl(\frac{h_2}{q\tilde{c}\tilde{r}f}\Bigr)\Bigr|. 
  \end{displaymath}
\end{lemma}

\textbf{Proof.} This is a direct consequence of the Poisson summation formula and Lemma \ref{char2}. \\

Finally we recall Weil's bound for Kloosterman sums
\begin{equation}\label{weil}
 | S(a, b, c)| \leq (a, b, c)^{1/2} c^{1/2} \tau( c). 
\end{equation} 
%as well as the trivial bound 
%\end{displaymath}
%   $|S(a, b, c)| \leq c$.   
By  twisted multiplicativity we see
\begin{equation}\label{twisted}
  S(qa, b, qc) = - S(a, b\bar{q}, c).
\end{equation}   
whenever $q \nmid bc$.  %We will later use the resulting  bounds
%   \begin{equation}\label{weil1}
  %  | S(a, b, qc)| \leq 2\sqrt{q} c, \quad q \nmid abc  \quad\text{and} \quad  |S(qa, b, qc) | \leq c, \quad q \nmid bc.
%   \end{equation}
   
%whenever $q \nmid abc$.  
\section{The main transformation}\label{main}
 
In this section we use Watson's formula and the Kuznetsov formula to transform the quantity of interest, $\sum_f \| f \|^4_4$, into character sums.
 
Let $f \in \mathcal{B}_q$ be an $L^2$-normalized cuspidal Hecke-Maa{\ss} newform of level $q$ with spectral parameter $t_f \leq T$. 
We begin with Parseval's identity 
\begin{displaymath}
 \| f \|_4^4 = \langle |f|^2, |f|^2 \rangle = V(q)^{-1} |\langle |f|^2, 1\rangle|^2 +  \sum_{g \in \mathcal{B}}| \langle |f|^2, g \rangle |^2+  \sum_{\mathfrak{a}} \int_{\Bbb{R}} |\langle |f|^2, E_{\mathfrak{a}}(., 1/2 + it) \rangle|^2 \frac{dt}{4\pi}. 
\end{displaymath}
We study the various terms on the right hand side. The constant function   contributes $V(q)^{-1} = O(1/q)$. Since the Laplace operator is symmetric, we have  
%\begin{displaymath}
  $\langle |f|^2, g \rangle =  (1/4 + t_g^2)^{-1} \langle \Delta |f|^2, g \rangle$. 
%\end{displaymath}
Iterating this procedure, we find $\langle |f|^2, g\rangle \ll_{T, A} t_g^{-A}$ for any $A \geq 0$. 
In particular, the oldforms  contribute
\begin{displaymath}
\begin{split}
  \sum_{g \in \mathcal{B} \setminus \mathcal{B}_q} |\langle |f|^2, g \rangle |^2 & = \sum_{\substack{ g \in \mathcal{B} \setminus \mathcal{B}_q\\ t_{g} \leq q^{\varepsilon}}} |\langle |f|^2, g \rangle |^2 + \sum_{\substack{ g \in \mathcal{B} \setminus \mathcal{B}_q\\ t_{g} > q^{\varepsilon}}} |\langle |f|^2, g \rangle |^2  \\
  & \ll_{\varepsilon, T}   \sum_{\substack{g \in \mathcal{B} \setminus \mathcal{B}_q\\ t_{g} \leq q^{\varepsilon}}} \| g \|_{\infty}^2  + \sum_{\substack{g \in \mathcal{B} \setminus \mathcal{B}_q\\ t_{g} > q^{\varepsilon}}}  (1+|t_{g}|)^{-3- \frac{1}{\varepsilon}} \ll q^{\varepsilon-1}
  \end{split}
\end{displaymath} 
by \eqref{weyl} and \eqref{upper}. By Watson's formula \cite[Theorem 5.1]{Wa} and positivity, the newforms contribute
\begin{displaymath}
\begin{split}
 \ll \frac{1}{q^2} \sum_{\substack{g \in \mathcal{B}_q\\ g \text{ even}}} \frac{\Lambda(1/2, f \times \bar{f} \times g)}{\Lambda(1, \text{Ad}^2 f)^2 \Lambda(1, \text{Ad}^2g)} & \ll_{\varepsilon, T}  \frac{1}{q^{2-\varepsilon}} \sum_{\substack{g \in \mathcal{B}_q\\ g \text{ even}}} \frac{L(1/2, f \times \bar{f} \times g)}{L(1, \text{Ad}^2 f) L(1, \text{Ad}^2 g)}e^{-\frac{3}{2} \pi |t_{g}|}.\\
 %& \leq \frac{1}{q^{2-\varepsilon}} \sum_{g \in \mathcal{B} } \frac{L(1/2, f \times f \times g)}{L(1, \text{Ad}^2 g)} e^{-\frac{3}{2} \pi t_{g}}. 
 \end{split}
\end{displaymath} 
Here we used a lower bound \cite{HL} on $L(1, \text{Ad}^2 f)$.  Since $f$ is an eigenfunction of the Fricke involution which is the scaling matrix for the cusp $\mathfrak{a} = 0$, the contribution of the two cusps is the same. By the unfolding technique we find as in \eqref{unfolding}
\begin{displaymath}
  \langle |f|^2, E_{\infty}(., s) \rangle = |\rho_f(1)|^2 \frac{2 L(s, f \times \bar{f})}{\zeta(2s)} \frac{  \Gamma_{\Bbb{R}}(s) \Gamma_{\Bbb{R}}(s - 2 i t_f) \Gamma_{\Bbb{R}}(s + 2 i t_f)}{2^{2+s}\Gamma_{\Bbb{R}}(1+s)}. 
\end{displaymath}
From \eqref{rho} we conclude that the Eisenstein contribution is 
\begin{displaymath}
  \ll_{\varepsilon, T}  \frac{1}{q^{2-\varepsilon}} \int_{\Bbb{R}} \frac{|L(1/2+it, f \times \bar{f})|^2}{|\zeta(1+2 i t)|^2} e^{-\frac{3}{2} \pi |t|} dt \ll q^{\varepsilon -1}
 \end{displaymath}
 by the convexity bound for $L(1/2+it, f \times \bar{f})$. Combining these estimates we find
 \begin{displaymath}
   \| f \|_4^4 \ll \frac{1}{q^{2-\varepsilon}L(1, \text{Ad}^2 f)} \sum_{\substack{g \in \mathcal{B}_q\\ g \text{ even}}} \frac{L(1/2, f \times \bar{f} \times g)}{L(1, \text{Ad}^2 g)} \bigl|e^{-\frac{3}{2} \pi t_{g}}\bigr| + q^{\varepsilon -1}.
 \end{displaymath}
We insert artificially the factor $G_1(0, t_g)G_2(0, t_g, t_f)$ by positivity and \eqref{lowerG} and also change the weight function $e^{-(3/2) \pi |t|}$ to the function 
\begin{displaymath}
   h(t) :=  \cosh\Bigl(\frac{t}{2A_2}\Bigr)^{-3 \pi A_2} \prod_{\nu=0}^{A_2}\Bigl(t^2 + \Bigl(\frac{1}{2} + \nu\Bigr)^2\Bigr).
 \end{displaymath}
Note that this function is holomorphic in $|\Im t|  < \pi A_2$ and has zeros at the zeros of $\cosh(\pi t)$ in this region.  Moreover, $h(t) \gg \exp(-\frac{3}{2} \pi |t|)$ for $t \in \mathcal{T}$.  Using \eqref{appr1} - \eqref{appr2}, we write 
\begin{displaymath}
  \| f \|^4_4 \ll \frac{q^{\varepsilon} }{q' L(1, \text{Ad}^2f)} \Bigl( 1+\sum_{\substack{g \in \mathcal{B}_q\\ g \text{ even}}} \frac{2 h(t_g)  }{q' L(1, \text{Ad}^2 g)}(1- \lambda_g(q)q^{1/2})S \Bigr)
\end{displaymath}  
%with $S_1+ S_2 \geq 0 $ given by  
where 
\begin{displaymath}
  S =  \sum_{n, m} \frac{\lambda_g(n)\lambda_{\text{Ad}^2f\times g}(m)}{(nm)^{1/2}} V_1\Bigl(\frac{n}{q^{1/2}}; t_g\Bigr)   V_2\Bigl(\frac{m}{q^{2}}; t_g, t_f\Bigr)
  \end{displaymath}
%  and 
 % \begin{displaymath}
 %   S_2 =  -q^{ -3/2} \sum_{\substack{g \in \mathcal{B}_q\\ g \text{ even}}} \frac{ e^{-\frac{3}{2} \pi t_{g}}  }{L(1, \text{Ad}^2 g)} \sum_{n, m} \frac{\lambda_g(q)\lambda_g(n)\lambda_{\text{Ad}^2f\times g}(m)}{(nm)^{1/2}} V_1\left(\frac{n}{q^{1/2}}; t_g\right)   V_2\left(\frac{m}{q^{2}}; t_g, t_f\right).
% \end{displaymath}
and $q'$ was defined in \eqref{deftildeq}. Hence by \eqref{upper} and \eqref{weyl} 
\begin{equation}\label{quantity}
  \sum_{t_f \leq T} \| f \|_4^4 \ll_{T, \varepsilon} q^{\varepsilon} +  q^{\varepsilon} \sum_{f \in \mathcal{B}_q}  \frac{2h(t_f)}{q' L(1, \text{Ad}^2 f)}  \sum_{\substack{g \in \mathcal{B}_q\\ g \text{ even}}} \frac{2 h(t_g) }{q' L(1, \text{Ad}^2 g)}(1- \lambda_g(q)q^{1/2})S . 
  %\sum_{\substack{g \in \mathcal{B}_q\\ g \text{ even}}} \frac{ \cosh(t_g/A)^{-A}  }{L(1, \text{Ad}^2 g)} \sum_{n, m} \frac{\lambda_g(n)\lambda_{\text{Ad}^2f\times g}(m)}{(nm)^{1/2}} V_1\left(\frac{n}{q^{1/2}}; t_g\right)   V_2\left(\frac{m}{q^{2}}; t_g, t_f\right)
\end{equation}
(Here we used \eqref{weyl} and $L(1, \text{Ad}^2 f) \gg_{T} q^{-\varepsilon}$ for the first term on the right hand side.) It is convenient to remove the terms with $q \mid nm$ in $S$.  By the rapid decay of $V_1$ the terms $q \mid n$ are negligible. Combining \eqref{small}  and \eqref{coeff} with \eqref{weyl},   \eqref{ranksel1}, \eqref{ranksel2} and the rapid decay of $V_2$, we see by trivial estimates  that the contribution of the terms $q \mid m$ in $S$ contributes at most $O(q^{\varepsilon - 1/4})$ to \eqref{quantity}.  Hence by \eqref{coeff} we are left with estimating
\begin{displaymath}
   \Sigma(q, q) := \Sigma_1(q, q) - \Sigma_2(q, q),
\end{displaymath}
say, where
\begin{displaymath}
\begin{split}
\Sigma_1(q, q) := & \sum_{f \in \mathcal{B}_q}  \frac{2h(t_f)}{q' L(1, \text{Ad}^2 f)} \sum_{\substack{g \in \mathcal{B}_q\\ g \text{ even}}} \frac{ 2h(t_g)  }{q' L(1, \text{Ad}^2 g)}  \sum_{q \nmid abdk} \frac{\mu(d) }{a b^2 k d^{3/2}} \\
  &  \times\sum_{q \nmid nm} \frac{\lambda_f(k^2)  \lambda_g(n)\lambda_f(m^2)\lambda_g(a^2dm)}{(nm)^{1/2}} V_1\Bigl(\frac{n}{q^{1/2}}; t_g\Bigr)   V_2\Bigl(\frac{a^2b^4k^2d^3m}{q^{2}}; t_g, t_f\Bigr)
 % & =: 
  % \sum_{q\nmid nm} \frac{\lambda_g(n)\lambda_{\text{Ad}^2f\times g}(m)}{(nm)^{1/2}} V_1\left(\frac{n}{q^{1/2}}; t_g\right)   V_2\left(\frac{m}{q^{2}}; t_g, t_f\right).
 \end{split} 
 \end{displaymath}
 and 
 \begin{displaymath}
\begin{split}
\Sigma_2(q, q) := &q^{1/2}  \sum_{f \in \mathcal{B}_q}  \frac{2h(t_f)}{q' L(1, \text{Ad}^2 f)} \sum_{\substack{g \in \mathcal{B}_q\\ g \text{ even}}} \frac{ 2h(t_g)  }{q' L(1, \text{Ad}^2 g)}  \sum_{q \nmid abdk} \frac{\mu(d)}{a b^2 k d^{3/2}} \\
  &  \times\sum_{q \nmid nm} \frac{ \lambda_f(k^2)  \lambda_g(qn)\lambda_f(m^2)\lambda_g(a^2dm)}{(nm)^{1/2}} V_1\Bigl(\frac{n}{q^{1/2}}; t_g\Bigr)   V_2\Bigl(\frac{a^2b^4k^2d^3m}{q^{2}}; t_g, t_f\Bigr). 
 % & =: 
  % \sum_{q\nmid nm} \frac{\lambda_g(n)\lambda_{\text{Ad}^2f\times g}(m)}{(nm)^{1/2}} V_1\left(\frac{n}{q^{1/2}}; t_g\right)   V_2\left(\frac{m}{q^{2}}; t_g, t_f\right).
 \end{split} 
 \end{displaymath}
 We would like to  apply the Kuznetsov formula to the spectral sums over $f$ and $g$. More precisely, we use \eqref{kuz1} for the sum over $f$ and \eqref{kuz3} for the sum over even Maa{\ss} forms $g$; in $\Sigma_2$ we use \eqref{kuz2} for the $f$-sum and \eqref{kuz4} for the $g$-sum.  However, this requires some preparation, as the $f$ and $g$-sum run only over cuspidal newforms, and both the oldforms of level 1 as well as the Eisenstein series are missing.  Therefore we add and subtract artificially the missing terms and define, in analogy with $\Sigma(q, q) = \Sigma_1(q, q) - \Sigma_2(q, q)$,  8 other quantities $\Sigma(*, *) = \Sigma_1(*, *) - \Sigma_2(*, *)$ where $* \in \{q, 1, \mathcal{E}\}$ in an obvious way in order to complete the spectral side of the Kuznetsov formula for the $f$ and $g$ sum respectively. For instance, we write
\begin{displaymath}
 \begin{split}
&  \Sigma_1(1, \mathcal{E}) := \sum_{f \in \mathcal{B}_1}  \frac{2h(t_f)}{(q+1)L(1, \text{Ad}^2 f)}    \int_{\Bbb{R}} \frac{ h(t) }{q''|\zeta^{(q)}(1+2it)|^2}  \sum_{q \nmid abdk} \frac{\mu(d) }{a b^2 k d^{3/2}} \\
  &\times   \sum_{q \nmid nm} \frac{c_1(f, q) \lambda_f(k^2) \eta(n,t)\lambda_f(m^2)\eta(a^2dm, -t)}{(nm)^{1/2}} V_1\Bigl(\frac{n}{q^{1/2}}; t\Bigr)   V_2\Bigl(\frac{a^2b^4k^2d^3m}{q^{2}}; t, t_f\Bigr) \frac{dt}{\pi},
  % \sum_{q\nmid nm} \frac{\lambda_g(n)\lambda_{\text{Ad}^2f\times g}(m)}{(nm)^{1/2}} V_1\left(\frac{n}{q^{1/2}}; t_g\right)   V_2\left(\frac{m}{q^{2}}; t_g, t_f\right).
 \end{split} 
\end{displaymath}
and
\begin{displaymath}
 \begin{split}
&  \Sigma_2(1, \mathcal{E}) := \sum_{f \in \mathcal{B}_1}  \frac{2h(t_f)}{(q+1)L(1, \text{Ad}^2 f)}    \int_{\Bbb{R}} \frac{ h(t)  }{q^{3/2}|\zeta^{(q)}(1+2it)|^2}  \sum_{q \nmid abdk} \frac{\mu(d) }{a b^2 k d^{3/2}} \\
  &\times   \sum_{q \nmid nm} \frac{ c_2(f, q) \lambda_f(k^2)\eta(n,-t)\lambda^{\ast}_f(m^2)\eta(a^2dm,-t)}{(nm)^{1/2}} V_1\Bigl(\frac{n}{q^{1/2}}; t\Bigr)   V_2\Bigl(\frac{a^2b^4k^2d^3m}{q^{2}}; t, t_f\Bigr) \frac{dt}{\pi},
  % \sum_{q\nmid nm} \frac{\lambda_g(n)\lambda_{\text{Ad}^2f\times g}(m)}{(nm)^{1/2}} V_1\left(\frac{n}{q^{1/2}}; t_g\right)   V_2\left(\frac{m}{q^{2}}; t_g, t_f\right).
 \end{split} 
\end{displaymath}
and similarly for all other combinations. We now apply the Kuznetsov formula to the completed expressions $\displaystyle \sum_{(\ast, \ast)} \Sigma_1(\ast, \ast)$ and   $\displaystyle\sum_{(\ast, \ast)} \Sigma_2(\ast, \ast)$, obtaining % In each case we get 4 terms $M_1^{\alpha, \beta}$ and $M_2^{\alpha, \beta}$ with $\alpha, \beta \in \{1, 2\}$ on the Kloosterman side as follows:
\begin{equation}\label{56}
  \Sigma(q, q) = - \sum_{\substack{(\ast, \ast) \in \{q, 1, \mathcal{E}\}^2\\ (\ast, \ast) \not= (q, q)}}( \Sigma_1(\ast, \ast) - \Sigma_2(\ast, \ast))  +  \sum_{q \nmid abdknm} \frac{\mu(d) }{(nm)^{1/2} a b^2 k d^{3/2}}  \sum_{\alpha, \beta, \gamma \in \{1, 2\} } M_{\alpha}^{\beta, \gamma}
\end{equation}
where %with
%\begin{displaymath}
 % \mathcal{V}(x, y; t_1, t_2) :=  \end{displaymath}
%one has
\begin{displaymath}
\begin{split}
  M_1^{1, 1} & = \delta_{n, a^2dm} \delta_{k, m} \mathcal{W}^{1, 1}\Bigl(\frac{n}{q^{1/2}},  \frac{a^2b^4k^2d^3m}{q^2}\Bigr)\\
  & \quad \mathcal{W}^{1, 1}(x, y) = \int_0^{\infty} \int_0^{\infty} V_1\left( x; t_2\right)   V_2\left( y; t_2, t_1\right) h(t_1)h(t_2) \frac{d^{\ast}t_1\,  d^{\ast}t_2}{2\pi^4},\\
     M_1^{1, 2} & =   \delta_{k, m} \sum_{q \mid c} \frac{1}{c} S(n,  a^2dm, c)\mathcal{W}^{1, 2}\Bigl(\frac{n}{q^{1/2}},  \frac{a^2b^4k^2d^3m}{q^2};  \frac{\sqrt{a^2dmn}}{c}\Bigr),\\
     & \quad \mathcal{W}^{1, 2}(x, y; \eta) = \int_0^{\infty} \int_0^{\infty} V_1\left( x; t_2\right)   V_2\left( y; t_2, t_1\right) h(t_1)h(t_2)\mathcal{J}^{0}\left(\eta, t_2\right) \frac{d^{\ast}t_1\,  d^{\ast}t_2}{\pi^3},\\
   \end{split}
   \end{displaymath}
   \begin{displaymath}
   \begin{split}
      M_1^{2, 1} & =   \delta_{n, a^2dm} \sum_{q \mid r} \frac{1}{r} S(m^2,  k^2, r) \mathcal{W}^{2, 1}\Bigl(\frac{n}{q^{1/2}},  \frac{a^2b^4k^2d^3m}{q^2};  \frac{km}{r}\Bigr),\\
   & \quad  \mathcal{W}^{2, 1}(x, y; \xi) = \int_0^{\infty} \int_0^{\infty} V_1\left( x; t_2\right)   V_2\left( y; t_2, t_1\right) h(t_1)h(t_2)\mathcal{J}^{+}\left(\xi, t_1\right) \frac{d^{\ast}t_1\,  d^{\ast}t_2}{2\pi^3},\\
     M_1^{2, 2} & =    \sum_{q \mid c}\sum_{q \mid r} \frac{1}{cr} S(n,  a^2dm, c) S(m^2, k^2, r)\mathcal{W}^{2, 2}\Bigl(\frac{n}{q^{1/2}},  \frac{a^2b^4k^2d^3m}{q^2}; \frac{km}{r}, \frac{\sqrt{a^2dmn}}{c}\Bigr),\\
     & \quad \mathcal{W}^{2, 2}(x, y; \xi, \eta) = \int_0^{\infty} \int_0^{\infty}  V_1\left( x; t_2\right)   V_2\left( y; t_2, t_1\right) h(t_1)h(t_2)\mathcal{J}^{+}\left(\xi, t_1\right)  \mathcal{J}^{0}\left(\eta, t_2\right) \frac{d^{\ast}t_1\,  d^{\ast}t_2}{\pi^2}.
  \end{split}
\end{displaymath}
Similarly,
\begin{displaymath}
\begin{split}
  M_2^{1, 1} &  = M_2^{2, 1} = 0, \\
   M_2^{1, 2} & =  q^{1/2}  \delta_{k, m} \sum_{q \mid c} \frac{1}{c} S(qn,  a^2dm, c)\mathcal{W}^{1, 2}\Bigl(\frac{n}{q^{1/2}},  \frac{a^2b^4k^2d^3m}{q^2};  \frac{\sqrt{a^2dmnq}}{c}\Bigr),\\
    M_2^{2, 2} & =   q^{1/2} \sum_{q \mid c}\sum_{q \mid r} \frac{1}{cr} S(qn,  a^2dm, c) S(m^2, k^2, r)\mathcal{W}^{2, 2}\Bigl(\frac{n}{q^{1/2}},  \frac{a^2b^4k^2d^3m}{q^2}; \frac{km}{r},  \frac{\sqrt{a^2dmnq}}{c}\Bigr).
       \end{split}
\end{displaymath}
In the rest of the paper we show that the $8+8+6 = 22$ (potentially) non-vanishing terms on the right hand side of \eqref{56} are all $O(q^{\varepsilon})$. This will complete the proof. %This is rather straightforward for the oldforms and Eisenstein terms $\Sigma_{1, 2}(\ast, \ast)$ and the diagonal or partially diagonal terms. The main work will be devoted to the off-off-diagonal terms $M_1^{2, 2}$ and $M_2^{2, 2}$. 

\section{The contribution of the oldforms and Eisenstein series}\label{old}

This section is devoted to bounding the terms $\Sigma_1(\ast, \ast)$ and $\Sigma_2(\ast, \ast)$ on the right hand side of  \eqref{56}.  All terms with $(\ast, \ast)\in \{1, \mathcal{E}\}^2$ can easily be bounded trivially: using only the Rankin-Selberg bounds \eqref{ranksel1}, \eqref{ranksel2} and the rapid decay of $V_1$ and $V_2$, we deduce
\begin{displaymath}
  \sum_{(\ast, \ast) \in \{1, \mathcal{E}\}^2} |\Sigma_1(\ast, \ast)| + |\Sigma_2(\ast, \ast)| \ll q^{-3/4 + \varepsilon}. 
\end{displaymath}

We proceed to bound the remaining terms $\Sigma_{1, 2}(q, \ast)$ and $\Sigma_{1, 2}(\ast, q)$ for $\ast \not = q$. The method for all these terms is identical, and we show as a typical example the case $\Sigma_1(q, \mathcal{E})$. By an inverse Mellin transform we have
\begin{displaymath}
 \begin{split}
&  \Sigma_1(q, \mathcal{E}) = \sum_{f \in \mathcal{B}_q}  \frac{2h(t_f)}{(q+1)L(1, \text{Ad}^2 f)}    \int_{\Bbb{R}} \frac{h(t)  }{q''|\zeta^{(q)}(1+2it)|^2} \\
&   \int_{(2)}\int_{(2)} \prod_{\pm} \Bigl(L^{(q)}(1/2 + u\pm it, \text{Ad}^2 f \times g)   \zeta^{(q)}(1/2 + v \pm it)\Bigr)   \widehat{V}_1\left(v; t\right)   \widehat{V}_2\left(u; t, t_f\right) q^{2u+\frac{v}{2}} \frac{du\, dv}{(2\pi i)^2}    \frac{dt}{\pi}. 
  % \sum_{q\nmid nm} \frac{\lambda_g(n)\lambda_{\text{Ad}^2f\times g}(m)}{(nm)^{1/2}} V_1\left(\frac{n}{q^{1/2}}; t_g\right)   V_2\left(\frac{m}{q^{2}}; t_g, t_f\right).
 \end{split} 
\end{displaymath}
We shift both contours to $\Re u = \Re v = \varepsilon$ and use the convexity bound $L(s, \text{Ad}^2 f \times g) \ll q^{1+\varepsilon}$ in $\Re s > 1/2$ (note that the poles of the zeta-function at $v=1/2 \pm it$ do not contribute by \eqref{zero}).  This yields the desired bound $\Sigma_1(q, \mathcal{E}) \ll q^{\varepsilon}$. The other 3 terms require only notational changes.\\

%We summarize the discussion so far by stating the preliminary bound
%\begin{displaymath}
 % \sum_{t_f \ll T}   \| f \|_4^4 \ll_{T, \varepsilon} q^{\varepsilon} + q^{\varepsilon}\Bigl| \sum_{q \nmid abdknm} \frac{ \mu(d)}{(nm)^{1/2} a b^2 k d^{3/2}} ( M_{1}^{2, 2} + M_2^{2, 2})\Bigr|. 
%\end{displaymath}
%Thus we need to bound
%\begin{equation}\label{M1}
% \mathcal{M}_1 := \sum_{q \nmid abdknm}  \sum_{q \mid c}\sum_{q \mid r} \frac{\mu(d) S(n,  a^2dm, c) S(m^2, k^2, r)}{(nm)^{1/2} a b^2 k d^{3/2}cr} \mathcal{W}^{2, 2}\Bigl(\frac{n}{q^{1/2}},  \frac{a^2b^4k^2d^3m}{q^2}; \frac{\sqrt{a^2dmn}}{c},   \frac{km}{r}\Bigr)
%\end{equation}
%and  
%\begin{equation}\label{M2}
 % \mathcal{M}_2 := q^{1/2} \sum_{q \nmid abdknm} \sum_{q\mid c}\sum_{q \mid r} \frac{\mu(d)S(qn,  a^2dm, c) S(m^2, k^2, r) }{(nm)^{1/2} a b^2 k d^{3/2}cr} \mathcal{W}^{2, 2}\Bigl(\frac{n}{q^{1/2}},  \frac{a^2b^4k^2d^3m}{q^2}; \frac{\sqrt{a^2dmnq}}{c},   \frac{km}{r}\Bigr).
%\end{equation}

\section{The weight functions}

In this technical section we provide useful bounds for the weight functions $\mathcal{W}$ occurring in the definition of the quantities $M_{\alpha}^{\beta, \gamma}$.  We start by collecting standard bounds for Bessel functions. The power series expansion implies
\begin{equation}\label{bessel1}
 e^{-\pi|t|} J_{2it}(x), \, e^{-\pi |t|}  I_{2it}(x) \ll_{\Im t}  (1+|t|)^{-1/2 + 2\Im t} x^{- 2\Im t}, \quad x \leq 1, t \in \Bbb{C}. 
\end{equation}
The asymptotic expansion implies
\begin{equation}\label{bessel2}
\begin{split}
&  \mathcal{J}^+(x, t) =  \frac{1}{\sqrt{x}} e\Bigl(\frac{x}{2\pi}\Bigr) v_+(x) +  \frac{1}{\sqrt{x}} e\Bigl(-\frac{x}{2\pi}\Bigr) v_-(x), \quad t \in \Bbb{R}, x \geq (1+|t|)^3. \\
% &  \mathcal{J}^0(x, t) \ll e^{-x}, \quad t \in \Bbb{R}, x \geq (1+|t|)^3,
 \end{split}
\end{equation}
where $v_{\pm}$ are smooth functions (depending on $t$) that satisfy $v_{\pm}^{(j)}(x) \ll_j  x^{-j}$ uniformly in $t$. We have for $j \in \Bbb{N}_0$ the general uniform upper bounds
\begin{equation}\label{bessel3}
\begin{split}
&  \frac{\partial^j}{\partial x^j} \mathcal{J}^{+}(x, t) \ll_j (1+|t|)^2 (1+x^{-j})  x^{-1/2}, \quad x > 0, t \in \Bbb{R},\\
& \frac{\partial^j}{\partial x^j} \mathcal{J}^-(x, t) \ll_{j, \varepsilon}  (1+x^{-j}) \times \begin{cases} x^{-\varepsilon},  & x < 1 + 10|t|, t \in \Bbb{R}\\
e^{-x/2}, & x > 1 + 10|t|, t \in \Bbb{R}. 
\end{cases}  
\end{split}
\end{equation}
These bounds are not optimal, but suffice for our application. 

Our first simple result shows that $\mathcal{W}$ is rapidly decreasing near $\infty$ in the first two variables and rapidly decreasing near $0$ in the other variables. 

\begin{lemma}\label{testfunct} The following uniform bounds hold for fixed $i, j\in \Bbb{N}_0$:
\begin{displaymath}
\begin{split}
    \mathcal{W}^{1, 1} (x, y) & \ll (1+x)^{-A_1} (1+y)^{-A_1},\\
    \mathcal{W}^{1, 2}(x, y; \xi) & \ll (1+x)^{-A_1} (1+y)^{-A_1}\min(\xi^{-1/2},  \xi^{4A_2}), \\
 \mathcal{W}^{2, 1}(x, y; \eta) & \ll (1+x)^{-A_1} (1+y)^{-A_1}\min(\eta^{-1/2},  \eta^{4A_2}) ,\\
   \mathcal{W}^{2, 2}(x, y; \xi,\eta) & \ll (1+x)^{-A_1} (1+y)^{-A_1}\min(\xi^{-1/2},  \xi^{4A_2})\min(\eta^{-1/2}, \eta^{4A_2}) . 
\end{split}  
\end{displaymath}
\end{lemma}

\textbf{Proof.} This follows directly by inserting the bounds from Lemma \ref{approx}. If $\xi$ and/or $\eta$ are greater than 1, we use \eqref{bessel3}; if $\xi$ and/or $\eta$ are less than 1, we write the corresponding  $t$-integral by symmetry as an integral over the whole real line,  shift the contour  down to $\Im t = - 2A_2$ (not crossing any poles) and use \eqref{bessel1}.  \\

We will also need the following more technical result. 
\begin{lemma}\label{refine} Let $N, M, Q, X \geq 1/2$,  and let $B  \in \Bbb{N}$, $\varepsilon > 0$ be fixed (but arbitrary). Let $\rho_1, \rho_2, \alpha_1, \alpha_2 > 0$ and $ \alpha_3 \in \Bbb{R}$ be real numbers and let $z, z_1, z_2 \in \Bbb{R}$. Let $w_1, w_2$ be two fixed smooth weight functions with support in $[1,2]$.  Then we have the uniform bounds
\begin{equation}\label{bound1}
\int_{\Bbb{R}} w_1\left(\frac{x}{N}\right) \mathcal{W}^{1, 2}(\rho_1 x, y; \sqrt{x} \alpha_1) e(-xz) dx \ll_B N(\sqrt{N}\alpha_1)^{-1/2} \Bigl(1 + Q^2\frac{|z|\sqrt{N}}{\alpha_1}\Bigr)^{-B}
\end{equation}
whenever $\alpha_1 \sqrt{N} \geq 1/Q$ and $y > 0$, and\footnote{Here the term $1/|\alpha_3|$ should be left out if $\alpha_3 = 0$, or one applies the convention $\min(x, \infty) = x$.}   
\begin{equation}\label{bound3}
\begin{split}
&\int_{\Bbb{R}} \int_{\Bbb{R}}  w_1\Bigl(\frac{x}{N}\Bigr)  w_2\Bigl(\frac{y}{M}\Bigr)\mathcal{W}^{2, 2}(\rho_1 x,  \rho_2 y;     \alpha_1y, \alpha_2\sqrt{yx})   e(-xz_1)  e(\alpha_3 y) e(-yz_2)dx\, dy\\
& \ll_{B, \varepsilon} \frac{XN^{1/2}}{\alpha_1^{1/2}\alpha_2}\Bigl(1+Q^2\frac{|z_1|\sqrt{N}}{\alpha_2\sqrt{M}}\Bigr)^{-B} \Bigl(1+Q^2|z_2|\min\Bigl(\frac{1}{\alpha_1}, \frac{1}{|\alpha_3|}, \frac{\sqrt{M}}{\alpha_2\sqrt{N}}\Bigr)\Bigr)^{-B} + MNX^{-B}
     \end{split}
     \end{equation}
whenever
\begin{equation}\label{assump}
 \min(\alpha_1 M, \alpha_2 \sqrt{NM} ) \geq \frac{1}{Q}, %\quad \max(\rho_1 N, \rho_2 M) \leq Q, 
 \quad  X\geq 10 + (\alpha_2\sqrt{NM})^{\varepsilon},
\end{equation}
and
 \begin{equation}\label{bound4}
\begin{split}
&\int_{\Bbb{R}}   w_2\Bigl(\frac{y}{M}\Bigr)   \mathcal{W}^{2, 2}( x,  \rho_2 y;     \alpha_1y, \alpha_2\sqrt{y})     e(-yz) dy \\
& \ll_{B, \varepsilon} \begin{cases} \frac{M^{1/4}}{(\alpha_1\alpha_2)^{1/2}} \Bigl(1+Q^2|z|\min\Bigl(\frac{1}{\alpha_1},  \frac{\sqrt{M}}{\alpha_2}\Bigr)\Bigr)^{-B}, & |z| \leq \alpha_1, \\
\frac{X}{\alpha_1^{1/2}\alpha_2} \Bigl(1+Q^2|z|\min\Bigl(\frac{1}{\alpha_1},  \frac{\sqrt{M}}{\alpha_2}\Bigr)\Bigr)^{-B} + MX^{-B}, & |z| \geq \alpha_1.\end{cases}
     \end{split}
     \end{equation}   
whenever
\begin{displaymath}
 \min(  \alpha_1 M,   \alpha_2 \sqrt{M} ) \geq \frac{1}{Q}, %\quad \max(\rho_1 N, \rho_2 M) \leq Q, 
 \quad  X\geq 10 +  (\alpha_2 \sqrt{M})^{\varepsilon}
\end{displaymath}
and $x > 0$. 
\end{lemma}

\textbf{Remark.} We will later apply this with $Q =  X = q^{\varepsilon}$, so as a first approximation  the reader can ignore the terms $Q^2$ and $MNX^{-B}$. \\

\textbf{Proof.} All three bounds depend on partial integration. We will always integrate the exponential factor containing $z, z_1, z_2$ respectively,  and differentiate all other factors. 

In order to prove \eqref{bound1}, we estimate trivially using \eqref{bessel3}, or we integrate by parts $B$ times and then estimate trivially using \eqref{bessel3}. Note that each integration by parts introduces an additional factor
\begin{displaymath}
  \frac{1}{|z|}\Bigl(\frac{1}{N}   + \frac{\alpha_1}{\sqrt{N}}\Bigr)\Bigl(1 + \frac{1}{\alpha_1\sqrt{N}}\Bigr) \leq \frac{1}{|z|} \Bigl((Q+ 1)\frac{\alpha_1}{\sqrt{N}}\Bigr)(1+Q) \ll Q^2 \frac{\alpha_1}{|z|\sqrt{N}}. 
\end{displaymath}\\

The  same strategy in the situation of \eqref{bound3} yields 
\begin{equation}\label{prelim}
 \frac{N^{3/4}M^{1/4}}{(\alpha_1\alpha_2)^{1/2}}   \Bigl(1+Q^2\frac{|z_1|\sqrt{N}}{\alpha_2 \sqrt{M}}\Bigr)^{-B}\Bigl(1+Q^2|z_2|\min\Bigl(\frac{1}{\alpha_1}, \frac{1}{|\alpha_3|}, \frac{\sqrt{M}}{\alpha_2\sqrt{N}}\Bigr)\Bigr)^{-B}.
  \end{equation}
 This bound suffices if $(\alpha_2 \sqrt{NM} )^{1/2} \leq X$. Let us now assume that
 \begin{displaymath}
   T := (\alpha_2 \sqrt{NM})^{1/2}/X\geq 1. 
 \end{displaymath}
 Then $T^{\varepsilon} \leq X^{1/2}$ by \eqref{assump}. If 
 \begin{equation}\label{assump1} 
    Q^2\frac{|z_1|\sqrt{N}}{\alpha_2 \sqrt{M}} \geq T^{\varepsilon/B} \quad \text{or} \quad Q^2|z_2|\min\Bigl(\frac{1}{\alpha_1}, \frac{1}{|\alpha_3|}, \frac{\sqrt{M}}{\alpha_2\sqrt{N}}\Bigr) \geq T^{\varepsilon/B},
 \end{equation}
 then we can replace $B$ by $B + B/\varepsilon$ in \eqref{prelim}, arriving at \eqref{bound3}. Let us now assume that \eqref{assump1} does not hold. Then
 \begin{displaymath}
    \Bigl(1+Q^2\frac{|z_1|\sqrt{N}}{\alpha_2 \sqrt{M}}\Bigr)^{-B}\Bigl(1+Q^2|z_2|\min\Bigl(\frac{1}{\alpha_1}, \frac{1}{|\alpha_3|},  \frac{\sqrt{M}}{\alpha_2\sqrt{N}}\Bigr)\Bigr)^{-B}\geq T^{-2\varepsilon}  \geq X^{-1},
 \end{displaymath}
 hence we only need to prove the upper bound $N^{1/2}/(\alpha_2\alpha_1^{1/2}) +MNX^{-B}$ for the double integral in \eqref{bound3}. Compared to the trivial estimate in \eqref{prelim} with $B=0$ we need to save a factor $(\alpha_2 \sqrt{NM})^{1/2}$. This comes from a standard stationary phase type argument. For convenience, we give precise details. We split the $t_2$-integral in the definition of $\mathcal{W}^{2, 2}$ into two pieces: $|t_2| \leq X^{2/3}$ and $|t_2| \geq X^{2/3}$. For large $t_2$, we estimate the $t_1, t_2$-integrals, as well as the above $x, y$-integral trivially using the rapid decay of the weight function $h$. This contributes the second term on the right hand side of \eqref{bound3}. For small $t_2$, we split $\mathcal{J}^0 = \frac{1}{2}(\mathcal{J}^+  +\mathcal{J}^{-})$. By \eqref{bessel3} we can bound the second term trivially due to the exponential decay of the Bessel-$K$-function getting again a contribution that is easily majorized by $\ll_{B} MNX^{-B}$. For $\mathcal{J}^+$ we  insert the asymptotic formula \eqref{bessel2}.  The $x$-integral  then becomes (for $y \asymp M$)
\begin{displaymath}
\begin{split}
 & \int_0^{\infty}  w_1\Bigl(\frac{x}{N}\Bigr) V_1(\rho_1 x; t_2) \frac{v_{\pm}(\alpha_2 \sqrt{xy}) }{\sqrt{\alpha_2 \sqrt{x y}}} e\Bigl(\pm \frac{\alpha_2\sqrt{xy}}{2\pi} - xz_1\Bigr)  dx =  \frac{N^{3/4}}{ \sqrt{\alpha_2\sqrt{y}}} \int_0^{\infty} W(x) e(\pm \beta_1 \sqrt{x} - \beta_2x) dx,   \end{split} 
\end{displaymath} 
say, where $W(x) = w_1(x) V_1(\rho_1N x; t_2) v_{\pm}(\alpha_2\sqrt{x N y}) x^{-1/2}$ is a function with support on $[1, 2]$ and bounded derivatives (uniformly in all parameters except $t_2$), and $\beta_1 = \alpha_2 \sqrt{Ny}/(2\pi)$, $\beta_ 2 = Nz_1$. If $|\beta_2/\beta_1|  \not\in [10^{-3}, 10^3]$, we integrate by parts sufficiently often, each time saving at least a factor $\alpha_2 \sqrt{Ny} \gg X^2$, and we obtain the trivial bound $O(MNX^{-B})$. If $|\beta_1| \asymp |\beta_2|$, then another change of variables yields
\begin{displaymath}
   \frac{N^{3/4}}{ \sqrt{\alpha_2\sqrt{y}}} \frac{\beta_1^2}{\beta_2^2} \int_0^{\infty} W\Bigl(x\frac{\beta_1^2}{\beta_2^2}\Bigr) e\Bigl(\frac{\beta_1^2}{|\beta_2|} (\pm \sqrt{x} - \text{sgn}(\beta_2) x)\Bigr) dx \ll  \frac{\sqrt{N}}{ \sqrt{\alpha_2\beta_1 \sqrt{y} }}  \asymp \frac{N^{1/2}}{\alpha_2M^{1/2}}
\end{displaymath}
by a standard stationary phase argument (e.g. \cite[p.\ 334]{St}). Integrating trivially over $y$ produces another factor $(M/\alpha_1)^{1/2}$, and the proof of \eqref{bound3} is complete in all cases. \\

The proof of \eqref{bound4} is almost identical, so we highlight only the key points. Integrating by parts and using \eqref{bessel3} yields a preliminary bound
\begin{displaymath}
  \frac{M^{1/4}}{(\alpha_1\alpha_2)^{1/2}}\Bigl(1+Q^2|z|\min(\frac{1}{\alpha_1}, \frac{\sqrt{M}}{\alpha_2}\Bigr)\Bigr)^{-B}. 
\end{displaymath}
This is acceptable if $X \geq (\alpha_2\sqrt{M})^{1/2}$ or if $|z| \leq \alpha_1$. In the other case, we argue as above, and hence we only need to show the upper bound $(\alpha_1^{1/2}\alpha_2)^{-1} + MX^{-B} $ for the integral in \eqref{bound4}. We cut the $t_1, t_2$-integral in the definition of $\mathcal{W}^{2, 2}$ according to whether $|t_1|$ and/or $|t_2|$ are bigger or less than $X^{2/3}$. By the rapid decay of the test function $h$, we can assume that both $t_1, t_2$ are small. In this range we can also replace $\mathcal{J}^0$ by $\mathcal{J}^+$ because of the rapid decay of the Bessel-$K$-function. For the two function $\mathcal{J}^+$ we insert the asymptotic expansion \eqref{bessel2}, and are left with the $y$-integral
\begin{displaymath}
  \int_0^{\infty} w_2\Bigl(\frac{y}{M}\Bigr) V_1(x; t_2) V_2(\rho_1 y; t_1, t_2)  \frac{v_{\pm}(\alpha_1y)}{\sqrt{\alpha_1 y}} \frac{v_{\pm}(\alpha_2\sqrt{y})}{\sqrt{\alpha_2\sqrt{y}}} e\Bigl(\pm \frac{\alpha_1 y}{2\pi}  \pm \frac{\alpha_2 \sqrt{y}}{2\pi} - yz\Bigr) dy.  
 \end{displaymath}
 By our present assumption $|z| \geq \alpha_1$ there is no phase cancellation in $(z\pm \alpha_1/(2\pi)) y$, and the same stationary phase argument for the $y$-integral followed by trivial estimates in the other integrals   gives  as before the bound \eqref{bound4}. 
   
% The proof of \eqref{bound3} is essentially identical and requires only notational changes. 
% x^{i_1}y^{i_2} \partial_x^{i_1} \partial_y^{i_2} \mathcal{W}^{1, 1} (x, y) & \ll (1+x)^{-A_1} (1+y)^{-A_1},\\
%   x^{i_1} y^{i_2} \partial_x^{i_1} \partial_y^{i_2} \partial^j_{\eta}  \mathcal{W}^{1, 2}(x, y, \eta) & \ll (1+x)^{-A_1} (1+y)^{-A_1}\min(\eta^{-1/2}, \eta^{4A_2-j}), \\
%  x^{i_1} y^{i_2} \partial_x^{i_1} \partial_y^{i_2}  \partial^j_{\xi} \mathcal{W}^{2, 1}(x, y, \xi) & \ll (1+x)^{-A_1} (1+y)^{-A_1}\min(\xi^{-1/2}, \eta^{4A_2-j}) ,\\
%   x^{i_1} y^{i_2} \partial_x^{i_1} \partial_y^{i_2} \partial^{j_1}_{\xi}\partial^{j_2}_{\eta} \mathcal{W}^{2, 2}(x, y, \xi,\eta) & \ll (1+x)^{-A_1} (1+y)^{-A_1}\min(\xi^{-1/2}, \xi^{4A_2-j_1})\min(\eta^{-1/2}, \eta^{4A_2-j_2}) . 

\section{Estimating character sums}

The scene has now been prepared to estimate the 6 potentially non-vanishing terms $M_{\alpha}^{\beta, \gamma}$ on the right hand side  of \eqref{56}. This is the heart of the proof of Theorem 1 and the most technical part. 

The bound
\begin{displaymath}
\sum_{q \nmid abdknm} \frac{1 }{(nm)^{1/2} a b^2 k d^{3/2}}( | M_1^{1, 1}| + |M_1^{1, 2}| + | M_1^{2, 1}|  ) \ll q^{\varepsilon}
\end{displaymath}
follows easily by trivial estimations using \eqref{weil} and the bounds from Lemma \ref{testfunct}. The other three terms need a more careful reasoning.

\subsection{The term $M_2^{1, 2}$}

Recall that we need to estimate
 % needs a more refined reasoning:
\begin{displaymath}
\begin{split}
&\sum_{q \nmid abdnm} \frac{\mu(d) M_2^{1, 2}}{n^{1/2} a b^2 m^{3/2} d^{3/2}}  = %\ll 1 +  q^{1/2} \Bigl|   
q^{1/2} \sum_{q \nmid abmdn} \sum_{q\mid c}
%& \times  \Bigl| \sum_{n \leq q^{1/2+\varepsilon}} \sum_{q \nmid a^2b^4m^3d^3 \leq q^{2+\varepsilon}} 
%\sum_{\substack{q \mid c\\ c \leq (a^2dmnq)^{1/2+\varepsilon}}} 
\frac{\mu(d) S(qn, a^2dm, c)}{n^{1/2} a b^2 m^{3/2} d^{3/2}c} \mathcal{W}^{1, 2}\Bigl(\frac{n}{q^{1/2}},  \frac{a^2b^4d^3m^3}{q^2};  \frac{a\sqrt{dmnq}}{c}\Bigr). 
\end{split}
\end{displaymath}
By the decay properties of $\mathcal{W}^{1, 2}$  given in Lemma \ref{testfunct} we can assume $a\sqrt{dmnq} \leq q^{7/4+\varepsilon}$, hence $q^2 \nmid c$.  Replacing $c$ by $cq$ with $q \nmid c$, we obtain by \eqref{twisted} (up to a negligible error)
\begin{displaymath}
- q^{-1/2} \sum_{q \nmid abmdnc} \frac{\mu(d)S(n, a^2dm\bar{q}, c) }{n^{1/2} a b^2 m^{3/2} d^{3/2}c} \mathcal{W}^{1, 2}\Bigl(\frac{n}{q^{1/2}},  \frac{a^2b^4d^3m^3}{q^2};  \frac{a\sqrt{dmn}}{cq^{1/2}}\Bigr).
\end{displaymath}
A trivial estimate gives only $O(q^{1/4+\varepsilon})$. In order to improve this, we can apply Poisson summation either in $a$ or in $n$, the latter being slightly  easier.  We can add   the terms $q \mid n$ with a negligible error, and we insert a smooth weight $w_1(n/N)w_2(a/A) w_3(c/C) (n/N)^{1/2} $ (using a smooth partition of unity) that localizes $N \leq n \leq 2N$, $A \leq a \leq 2A$ and $C \leq c \leq 2C$. Again by Lemma \ref{testfunct} we can assume
\begin{equation}\label{sizes}
  %N \leq q^{1/2+\varepsilon}, \quad 
  A \leq \frac{q^{1+\varepsilon}}{b^2 (dm)^{3/2}}, \quad C \leq \frac{A\sqrt{dmN}}{q^{1/2-\varepsilon}}.
\end{equation}
 Thus we need to bound
\begin{equation}\label{needtobound}
\begin{split}
&  \frac{1}{q^{1/2}N^{1/2}} \sum_{q \nmid abmdc} \frac{w_2(a/A)w_3(c/C) }{a  b^2 m^{3/2} d^{3/2}c} \\
  & \times \Bigl|\sum_{n} S(n, a^2dm\bar{q}, c) w_1\left(\frac{n}{N}\right) \mathcal{W}^{1, 2}\Bigl(\frac{n}{q^{1/2}},  \frac{a^2b^4d^3m^3}{q^2};  \frac{a\sqrt{dmn}}{cq^{1/2}}\Bigr)\Bigr|.
\end{split}  
\end{equation}
By Lemma \ref{char1} with $\alpha = a^2 dm\bar{q}$ and  $\gamma = c$, the $n$-sum is
% By Poisson summation, the inner sum equals
%\begin{equation}\label{nsum}
  %\frac{1}{c} \sum_{\nu \, ( c)} S(\nu, a^2 dm\bar{q}, c) \sum_{h \in \Bbb{Z}} e\left(\frac{h\nu}{c}\right) F\left(\frac{h}{c}\right)  = \sum_{h\not=0} e\left(-\frac{a^2dm\bar{h}}{c}\right) F\left(\frac{h}{c}\right) \leq \sum_{h \not= 0} \left| F\left(\frac{h}{c}\right) \right|
%\end{equation}
%where
\begin{equation}\label{nsum}
 \leq \sum_{h \not= 0} \Bigl|   \int_0^{\infty} w_1\left(\frac{x}{N}\right)  \mathcal{W}^{1, 2}\Bigl(\frac{x}{q^{1/2}},  \frac{a^2b^4d^3m^3}{q^2};  \frac{a\sqrt{dmx}}{cq^{1/2}}\Bigr) e\Bigl(-\frac{xh}{c}\Bigr) dx\Bigr|. 
\end{equation}
The estimate \eqref{bound1} with
\begin{displaymath}
  \rho_1 = \frac{1}{q^{1/2}}, \quad \alpha_1 = \frac{a\sqrt{dm}}{cq^{1/2}}, \quad Q = q^{\varepsilon}
\end{displaymath}
(and $y = a^2b^4d^3m^3/q^2$) is applicable by \eqref{sizes} if $a\asymp A$ and $c \asymp C$ and  %implies %Integrating by parts sufficiently often using Lemma \ref{testfunct}, we have
bounds \eqref{nsum} by
\begin{displaymath}
\begin{split}
 &  q^{\varepsilon}\sum_{h \not = 0} N \Bigl(\frac{a\sqrt{dmN}}{cq^{1/2}}\Bigr)^{-1/2}  \Bigl(1 + |h| \sqrt{N} \frac{q^{1/2}}{a\sqrt{dm}}\Bigr)^{-10}\\
%\end{displaymath} 
%Thus \eqref{nsum} is at most
%\begin{displaymath}
 & \ll q^{\varepsilon} N \Bigl(\frac{a\sqrt{dmN}}{cq^{1/2}}\Bigr)^{-1/2} \Bigl(\sqrt{N} \frac{q^{1/2}}{a\sqrt{dm}}\Bigr)^{-1} = \frac{N^{1/4} c^{1/2}a^{1/2}(dm)^{1/4} }{q^{1/4-\varepsilon}}.
\end{split} 
\end{displaymath}
We substitute this back into \eqref{needtobound} getting the final bound 
\begin{displaymath}
\frac{1}{q^{3/4-\varepsilon}N^{1/4}}  \sum_{  abmdc} \frac{w_2(a/A)w_3(c/C) }{(ac)^{1/2}  b^2 (md)^{5/4}}  \ll q^{\varepsilon}
\end{displaymath}
by \eqref{sizes}. 
%Using \eqref{sizes}, it is easy to bound \eqref{needtobound} by $O(q^{\varepsilon})$. \\

%The remaining sections are dedicated to bound the two off-off-diagonal terms. 

\subsection{The   term $M_1^{2, 2}$}

%By the  decay properties of $\mathcal{W}^{2, 2}$   we can assume $c \leq q^{5/4+\varepsilon}$, hence $q^2 \nmid c$, and $r \leq q^{2+\varepsilon}$. We estimate trivially the terms with $q^2 \mid r$ by \eqref{weil}:
%\begin{displaymath}
%  \sum_{\substack{q \nmid n\\ q \leq q^{1/2+\varepsilon}}} \sum_{\substack{q \nmid abkdm\\ a^2b^4k^2d^3m \ll q^{2+\varepsilon}}} \sum_{\substack{q \mid c\\ c \leq a\sqrt{dmn}q^{\varepsilon}}} \sum_{\substack{q^2 \mid r\\ r \leq kmq^{\varepsilon}}} 
 % \sum_{q \nmid abdknm}  \sum_{q \mid c}\sum_{q^2 \mid r} \frac{(n,  a^2dm, c) (m^2, k^2, r)}{(nm)^{1/2} a b^2 k d^{3/2}(cr)^{1/2}} \ll q^{-1/8+\varepsilon}.
%\end{displaymath}
%Now we replace $c$ with $qc$ and $r$ with $qr$ where $q \nmid cr$ and obtain
%\begin{displaymath}
%q^{-2} \sum_{q \nmid abdknmcr}   \frac{\mu(d) S(n,  a^2dm, qc) S(m^2, k^2, qr)}{(nm)^{1/2} a b^2 k d^{3/2}cr} \mathcal{W}^{2, 2}\Bigl(\frac{n}{q^{1/2}},  \frac{a^2b^4k^2d^3m}{q^2}; \frac{\sqrt{a^2dmn}}{qc},   \frac{km}{qr}\Bigr). 
%\end{displaymath}

Here we need to bound
\begin{displaymath}
\frac{1}{q^2}  \sum_{q \nmid abdknm}  \sum_{c, r} \frac{\mu(d) S(n,  a^2dm, qc) S(m^2, k^2, qr)}{(nm)^{1/2} a b^2 k d^{3/2}cr} \mathcal{W}^{2, 2}\Bigl(\frac{n}{q^{1/2}},  \frac{a^2b^4k^2d^3m}{q^2};  \frac{km}{qr}, \frac{\sqrt{a^2dmn}}{qc}\Bigr). 
\end{displaymath}
 The key variables  are then $n, m, c, r$, and the reader can savely think of the other variables as 1. We re-include the terms $q \mid m$. By the decay properties of $\mathcal{W}^{2, 2}$,  \eqref{weil}, \eqref{twisted} and trivial estimates, these contribute $O(q^{-9/8+\varepsilon})$. We can also include the terms $q \mid n$ at a negligible cost. It is convenient to include smooth weights 
\begin{displaymath}
  w_1(n/N)w_2(m/M) w_3(c/C) w_4(r/R) \frac{(nm)^{1/2}cr}{(NM)^{1/2} CR}   
  \end{displaymath} 
   where all $w_j$ have support on $[1,2]$, and the parameters $N, M, R, C \geq 1/2$ satisfy (cf.\  Lemma \ref{testfunct})
\begin{equation}\label{sizes1}
  N  \leq q^{1/2+\varepsilon}, \quad M \leq \frac{q^{2+\varepsilon}}{a^2b^4k^2d^3},  \quad R \leq \frac{kM}{q^{1-\varepsilon}}, \quad C \leq \frac{a\sqrt{dMN}}{q^{1-\varepsilon}}. 
\end{equation}
The idea is now to apply Poisson summation in the $n$ and $m$ variable. Since the $m$-variable is very long, the second application is certainly advantageous. The benefit of the first application is not immediately obvious, since it makes the $n$-sum (which is generically of length $q^{1/2}$) longer (the new $h_1$-sum is generically of length $q^{3/4}$). The point here is that the $n$-sum is a linear exponential sum, and hence the resulting complete double sum after both applications of Poisson simplifies a lot which compensates the loss in the length. 

More formally, we now apply Lemma \ref{char1} with $\alpha = a^2 d m$ and $\gamma = qc$ to the $n$-sum, and then apply Lemma \ref{char3} with $\beta = -a^2d\bar{h}_1$ and $\kappa = k^2$ to the $m$-sum.  Unfortunately this introduces a zoo of new variables. As in Lemma  \ref{char3} we decompose $r = f\tilde{r}, c = f\tilde{c}$ with $f = (r, c)$. Moreover, we decompose $f = f_1f_2^2$, $\tilde{r} = \tilde{r}_1 \tilde{r}_2^2$ with $\mu(f_1)^2 = \mu(\tilde{r}_1)^2 = 1$, so that 
\begin{equation}\label{para}
\begin{split}
&c = \tilde{c}f_1f_2^2,\\
  &r = \tilde{r}_1\tilde{r}_2^2f_1f_2^2 = \frac{f_1\tilde{r}_1}{(f_1, \tilde{r}_1)^2} \times  (f_1, \tilde{r}_1)^2 f_2^2 \tilde{r}_2^2, \quad \mu\Bigl(\frac{ f_1\tilde{r}_1}{(f_1, \tilde{r}_1)^2}\Bigr)^2 = 1.
 \end{split} 
\end{equation}
 In this way we obtain the upper bound 
\begin{displaymath}
\begin{split}
 &   \frac{1}{q^2 } \sum_{q \nmid abdk} \sum_{f_1f_2\tilde{c}\tilde{r}_1\tilde{r}_2} \frac{ (f_1, \tilde{r}_1)f_2\tilde{r}_2}{a b^2 k d^{3/2} (NM)^{1/2} CR} w_3\Bigl(\frac{f_1f_2^2\tilde{c}}{C}\Bigr) w_4\Bigl(\frac{f_1f_2^2\tilde{r}_1\tilde{r}_2^2}{R}\Bigr)  \\ %\sum_{\substack{h_1 \in \Bbb{Z}\\ (h_1, qf_1f_2^2\tilde{c}) = 1}} %\sum_{\substack{h_2 \in \Bbb{Z}\\ h_1h_2 \equiv a^2d \tilde{r}_1\tilde{r}_2^2 \, (\tilde{c})}}
%\end{displaymath}
%\begin{displaymath}
%\times  
% \sum_{\substack{h_1 \in \Bbb{Z}\\ (h_1, qc) = 1}} 
&%\sum_{\substack{h_2 \in \Bbb{Z}\\ h_1h_2 \equiv a^2d \tilde{r}_1\tilde{r}_2^2 \, (\tilde{c})}}
\quad \times \underset{\substack{h_1,  h_2 \in \Bbb{Z}\\  (h_1, qf_1f_2^2\tilde{c}) = 1\\ h_1h_2 \equiv a^2d \tilde{r}_1\tilde{r}_2^2 \, (\tilde{c})}}{\sum\sum} \Bigl|  \int_{\Bbb{R}}\int_{\Bbb{R}}  \mathcal{W}^{2, 2}\Bigl(\frac{x}{q^{1/2}},  \frac{a^2b^4k^2d^3y}{q^2}; \frac{ky}{qf_1f_2^2\tilde{r}_1\tilde{r}_2^2}, \frac{\sqrt{a^2dxy}}{qf_1f_2^2\tilde{c}}  \Bigr)\\
 & \quad\quad  \quad\quad \times w_1\Bigl(\frac{x}{N}\Bigr)w_2\Bigl(\frac{y}{M}\Bigr)  e\Bigl(-\frac{xh_1}{qf_1f_2^2\tilde{c}}\Bigr) e\Bigl(-\frac{yh_2}{q\tilde{c} \tilde{r}_1\tilde{r}_2^2f_1f_2^2}\Bigr)  dx\, dy\Bigr|. 
 \end{split}
\end{displaymath}
The bound \eqref{bound3}  with 
\begin{displaymath}
\begin{split}
 & \rho_1 = \frac{1}{q^{1/2}}, \quad \rho_2 = \frac{a^2b^4k^2d^3}{q^2}, \quad \alpha_1 = \frac{k}{qf_1f_2^2\tilde{r}_1\tilde{r}_2^2}, \quad \alpha_2 = \frac{a \sqrt{d}}{qf_1f_2^2\tilde{c}}, \quad \alpha_3 = 0, \\
  & z_1 = \frac{h_1}{qf_1f_2^2 \tilde{c}}, \quad z_2 = \frac{h_2}{q\tilde{c} \tilde{r}_1 \tilde{r}_2^2 f_1f_2^2}, \quad Q = X = q^{\varepsilon}
\end{split}  
\end{displaymath}
is applicable by \eqref{sizes1} if $f_1f_2^2 \tilde{c} \asymp C$ and $f_1f_2^2 \tilde{r}_1\tilde{r}_2^2 \asymp R$ and implies that the double integral is at most
\begin{displaymath}
 \ll q^{\varepsilon} \frac{N^{1/2}q^{3/2} CR^{1/2}}{a (dk)^{1/2} } \Bigl(1+\frac{|h_1|\sqrt{N} }{a\sqrt{dM}}\Bigr)^{-10} \Bigl(1+ |h_2| \min\Bigl(\frac{1}{k \tilde{c}}, \frac{\sqrt{M}C}{a \tilde{c}R \sqrt{dN}} \Bigr)\Bigr)^{-10} + q^{-100}.  
  \end{displaymath}
We can ignore the second term, and we sum the first term  over $h_1, h_2$ getting the upper bound
\begin{displaymath}
  \ll q^{\varepsilon} \frac{N^{1/2}q^{3/2} CR^{1/2}}{a (dk)^{1/2} } \frac{a\sqrt{dM}}{\sqrt{N}}  \Bigl(k + \frac{aR\sqrt{dN}}{\sqrt{M}C}\Bigr).
\end{displaymath}
We substitute this back and sum over $f_1, f_2, \tilde{r}_1, \tilde{r}_2, \tilde{c}$ and obtain a total contribution of
\begin{displaymath}
\begin{split}
 & \frac{1}{q^{2-\varepsilon} } \sum_{  abdk}  \frac{1}{a b^2 k d^{3/2} (NM)^{1/2} } \frac{N^{1/2}q^{3/2} CR^{1/2}}{a (dk)^{1/2} } \frac{a\sqrt{dM}}{\sqrt{N}}  \Bigl(k + \frac{aR\sqrt{dN}}{\sqrt{M}C}\Bigr)\\
  & = \frac{1}{q^{1/2-\varepsilon} } \sum_{abdk} \Bigl(\frac{CR^{1/2}}{ab^2d^{3/2}k^{1/2}N^{1/2}}  + \frac{R^{3/2}}{b^2dk^{3/2}M^{1/2}}\Bigr).
 \end{split} 
\end{displaymath}
We insert the upper bound for $C$ and $R$ from \eqref{sizes1}, then the upper bound for $M$, and obtain the desired bound $q^{\varepsilon}$. 

\subsection{The   term $M_2^{2, 2}$}

We argue similarly as in the previous subsection and consider the term
\begin{displaymath}
 \frac{1}{q^{3/2}} \sum_{q \nmid abdknm} \sum_{c, r} \frac{\mu(d)S(qn,  a^2dm, qc) S(m^2, k^2, qr) }{(nm)^{1/2} a b^2 k d^{3/2}cr} \mathcal{W}^{2, 2}\Bigl(\frac{n}{q^{1/2}},  \frac{a^2b^4k^2d^3m}{q^2};  \frac{km}{qr}, \frac{a\sqrt{dmn}}{\sqrt{q}c}\Bigr).
\end{displaymath}
First we observe that by the decay properties of $\mathcal{W}^{2, 2}$  from Lemma \ref{testfunct}  we can assume that $a\sqrt{ d mn /q} \leq q^{3/4+\varepsilon}$, hence $q \nmid c$. %Again by the decay properties of $\mathcal{W}^{2, 2}$ we can assume $r \leq q^{2+\varepsilon}$, hence $q^3 \nmid r$. We estimate trivially the terms with $q^2 \mid r$ as $O()$. 
We rewrite the previous display using \eqref{twisted} (up to a negligible error and up to sign) as
\begin{displaymath}
   q^{-3/2} \sum_{q \nmid abdknmc} \sum_{r}   \frac{\mu(d)S(n,  a^2dm\bar{q}, c) S(m^2, k^2, qr) }{(nm)^{1/2} a b^2 k d^{3/2}cr} \mathcal{W}^{2, 2}\Bigl(\frac{n}{q^{1/2}},  \frac{a^2b^4k^2d^3m}{q^2};  \frac{km}{qr}, \frac{a\sqrt{dmn}}{\sqrt{q}c}\Bigr).
\end{displaymath}
Again the key players are the variables $n, m, r, c$. We re-introduce the terms $q\mid m$ which by \eqref{weil} and trivial estimates infers an error $O(q^{-3/8+\varepsilon})$. The terms $q \mid n$ can be included with a negligible error. Next we introduce smooth weights
\begin{displaymath}
  w_1(n/N)w_2(m/M) w_3(c/C) w_4(r/R) \frac{(nm)^{1/2}cr}{(NM)^{1/2} CR}   
  \end{displaymath} 
   where all $w_j$ have support on $[1,2]$, and the parameters $N, M, R, C \geq 1/2$ satisfy (cf.\  Lemma \ref{testfunct})
\begin{equation}\label{sizes2}
  N  \leq q^{1/2+\varepsilon}, \quad M \leq \frac{q^{2+\varepsilon}}{a^2b^4k^2b^3}, \quad R \leq \frac{kM}{q^{1-\varepsilon}}, \quad C \leq \frac{a\sqrt{dMN}}{q^{1/2-\varepsilon}}. 
\end{equation}
There is one special case that we need to treat separately: if $c = 1$, then the Kloosterman sum $S(n, a^2dm\bar{q}, c)$ degenerates. We will postpone this case for the moment and assume $C > 1$ so that automatically $c \not=1$. Now we apply Lemma \ref{char1} with $\gamma = c$, $\alpha= a^2 d m \bar{q}$, getting
 \begin{displaymath}
 \begin{split}
   \frac{1}{q^{3/2}} & \sum_{q \nmid abdkc} \sum_{m, r} \sum_{(h_1, c) = 1}   \frac{\mu(d) w_2(m/M)w_3(c/C)w_4(r/R)}{(NM)^{1/2} a b^2 k d^{3/2}CR}  e\Bigl(-\frac{a^2dm\overline{qh_1}}{c}\Bigr) S(m^2, k^2, qr) \\
   & \times \int_{\Bbb{R}} w_1\Bigl(\frac{x}{N}\Bigr) \mathcal{W}^{2, 2}\Bigl(\frac{x}{q^{1/2}},  \frac{a^2b^4k^2d^3m}{q^2};  \frac{km}{qr}, \frac{a\sqrt{dmx}}{\sqrt{q}c}\Bigr) e\left(-\frac{xh_1}{c}\right) dx.
 \end{split}  
\end{displaymath}
Since $c \not= 1$, we have $h_1 \not = 0$. Hence we can use the reciprocity formula
\begin{displaymath}
  e\Bigl(-\frac{a^2dm\overline{qh_1}}{c}\Bigr)  = e\Bigl(\frac{a^2dm\bar{c}}{qh_1}\Bigr) e\Bigl(-\frac{a^2dm}{cqh_1}\Bigr). 
\end{displaymath}
In order to display the similarities to the computation in the previous subsection, we switch the roles of $c$ and $h_1$, obtaining
 \begin{displaymath}
 \begin{split}
   \frac{1}{q^{3/2}} & \sum_{q \nmid abdk} \sum_{\substack{m, r\\ c \in \Bbb{Z}\setminus \{0\}}} \sum_{(h_1, qc) = 1}   \frac{\mu(d) w_2(m/M)w_3(h_1/C)w_4(r/R)}{(NM)^{1/2} a b^2 k d^{3/2}CR}  e\Bigl(\frac{a^2dm\overline{h_1}}{qc}\Bigr) S(m^2, k^2, qr) \\
   & \times e\Bigl(-\frac{a^2dm}{cqh_1}\Bigr) \int_{\Bbb{R}} \mathcal{W}^{2, 2}\Bigl(\frac{x}{q^{1/2}},  \frac{a^2b^4k^2d^3m}{q^2};     \frac{km}{qr}, \frac{a\sqrt{dmx}}{\sqrt{q}h_1}\Bigr) w_1\Bigl(\frac{x}{N}\Bigr) e\left(\frac{xc}{h_1}\right) dx.
 \end{split}  
\end{displaymath}
This having done, we now apply Lemma \ref{char3} with $\kappa = k^2$, $\beta = a^2d \bar{h}_1$ and use the same parametrization \eqref{para} as in the previous estimation. Thus we arrive at the upper bound
\begin{displaymath}
\begin{split}
 \frac{1}{q^{3/2}} &\sum_{q \nmid abdk} \sum_{\substack{\tilde{r}_1, \tilde{r}_2, f_1, f_2\\ \tilde{c} \in \Bbb{Z}\setminus \{0\}}} \frac{(f_1, \tilde{r}_1)f_2\tilde{r}_2}{(NM)^{1/2} a b^2 k d^{3/2}CR} \underset{\substack{h_1\in \Bbb{N}, h_2 \in \Bbb{Z}\\ (h_1, qf_1f_2^2\tilde{c}) = 1\\h_1h_2 \equiv a^2 d \tilde{r}_1\tilde{r}^2_2\, (\tilde{c})}}{\sum\sum}  
      w_3\Bigl(\frac{h_1}{C}\Bigr)w_4\Bigl(\frac{\tilde{r}_1\tilde{r}_2^2f_1f_2^2}{R}\Bigr)   \\
   & \times \Bigl| \int_{\Bbb{R}} \int_{\Bbb{R}} \mathcal{W}^{2, 2}\Bigl(\frac{x}{q^{1/2}},  \frac{a^2b^4k^2d^3y}{q^2}; \frac{ky}{qf_1f_2^2\tilde{r}_1\tilde{r}_2^2}, \frac{a\sqrt{dyx}}{\sqrt{q}h_1}\Bigr) w_1\Bigl(\frac{x}{N}\Bigr)  w_2\Bigl(\frac{y}{M}\Bigr)\\ 
   &\quad\quad\quad 
     \times e\left(-\frac{xf_1f_2^2 \tilde{c}}{h_1}\right)  e\Bigl(-\frac{a^2dy}{\tilde{c}f_1f_2^2qh_1}\Bigr) e\Bigl(-\frac{yh_2}{q\tilde{c}\tilde{r}_1\tilde{r}_2^2f_1f_2^2}\Bigr)dx\, dy\Bigr|.
\end{split}
\end{displaymath}
The bound \eqref{bound3} with
\begin{displaymath}
  \begin{split}
 & \rho_1 = \frac{1}{q^{1/2}}, \quad \rho_2 = \frac{a^2b^4k^2d^3}{q^2}, \quad  \alpha_1 = \frac{k}{qf_1f_2^2\tilde{r}_1\tilde{r}_2^2}, \quad \alpha_2 = \frac{a \sqrt{d}}{\sqrt{q}h_1},  \quad \alpha_3 = - \frac{a^2d}{\tilde{c} f_1f_2^2 q h_1},\\
  & z_1 = \frac{f_1f_2^2\tilde{c}}{h_1}, \quad z_2 = \frac{h_2}{q\tilde{c} \tilde{r}_1 \tilde{r}_2^2 f_1f_2^2}, \quad Q = X = q^{\varepsilon}
\end{split}  
\end{displaymath}
is applicable by \eqref{sizes2} and implies that  the double integral is at most
\begin{displaymath}
q^{\varepsilon} \frac{ N^{1/2}qCR^{1/2}}{ a (dk)^{1/2}} \Bigl(1 + f_1f_2^2|\tilde{c}| \frac{\sqrt{qN}}{a\sqrt{dM}}\Bigr)^{-10} \Bigl(1 + |h_2| \min\Bigl(\frac{1}{k|\tilde{c}|}, \frac{C\sqrt{M}}{a\sqrt{dNq} R |\tilde{c}|}, \frac{Cf_1f_2^2}{a^2dR}\Bigr)  \Bigr)^{-10},
 \end{displaymath}
up to a negligible term $q^{-100}$. Now it's just a matter of book-keeping. The sum over $h_2$ is at most
\begin{displaymath}
  \ll k + \frac{a\sqrt{dNq} R }{C\sqrt{M}} + \frac{a^2dR}{C|\tilde{c}|f_1f_2^2}. 
\end{displaymath}
We sum this over $\tilde{c}$ and then over $h_1$ and $\tilde{r}_1\tilde{r}_2^2f_1f_2^2$  getting
\begin{displaymath}
  CR \Bigl(\frac{a\sqrt{dM}}{\sqrt{qN}}\Bigl(k + \frac{a\sqrt{dNq} R }{C\sqrt{M}}\Bigr) + \frac{a^2dR}{C}\Bigr),
\end{displaymath}
so that the total contribution is given by
\begin{displaymath}
\begin{split}
&\frac{1}{q^{3/2-\varepsilon}} \sum_{adbk} \frac{N^{1/2}qCR^{1/2}}{(NM)^{1/2} (abd)^2k^{3/2}} \Bigl(\frac{a\sqrt{dM}}{\sqrt{qN}}\Bigl(k + \frac{a\sqrt{dNq} R }{C\sqrt{M}}\Bigr) + \frac{a^2dR}{C}\Bigr)\\
& = \sum_{adbk} \Bigl(\frac{CR^{1/2}}{ab^2d^{3/2}k^{1/2}qN^{1/2}} + \frac{R^{3/2}}{b^2dk^{3/2}(qM)^{1/2}} \Bigr).
\end{split}
\end{displaymath}
The desired bound $q^{\varepsilon}$ follows now easily from \eqref{sizes2}. \\

It remains to treat the case $c=1$, that is,
\begin{displaymath}
     \sum_{q \nmid abdk} \sum_{m, n, r}   \frac{w_1(n/N)w_2(m/M) w_4(r/R)}{q^{3/2} (MN)^{1/2}R}   \frac{\mu(d) S(m^2, k^2, qr) }{  a b^2 k d^{3/2}} \mathcal{W}^{2, 2}\Bigl(\frac{n}{q^{1/2}},  \frac{a^2b^4k^2d^3m}{q^2};  \frac{km}{qr}, \frac{a\sqrt{dmn}}{\sqrt{q}}\Bigr).
\end{displaymath}
where $M, N, R$ are subject to \eqref{sizes2}, and we may also assume $a \sqrt{dMN/q} \geq q^{-\varepsilon}$. The argument in this special case is not much different from the general case above, and a little easier.  We apply Lemma \ref{char3} with $c = 1$, $\beta = 0$ and $\kappa = k^2$. Writing $r = r_1r_2^2$ with $\mu(r_1)^2 = 1$, we obtain the upper bound
\begin{displaymath}
\begin{split}
  &   \sum_{q \nmid abdk} \sum_{ \substack{n, r_1, r_2\\ h \in \Bbb{Z}}}   \frac{w_1(n/N)  w_4(r_1r_2^2/R) r_2}{q^{3/2} (MN)^{1/2}R  a b^2 k d^{3/2}} \\
     & \times  \int_0^{\infty}  w_2\Bigl(\frac{y}{M}\Bigr) \mathcal{W}^{2, 2}\Bigl(\frac{n}{q^{1/2}},  \frac{a^2b^4k^2d^3y}{q^2};  \frac{ky}{qr_1r_2^2}, \frac{a\sqrt{dyn}}{\sqrt{q}}\Bigr) e\Bigl(-\frac{yh}{qr_1r_2^2}\Bigr) dy.
  \end{split}   
\end{displaymath}
We apply \eqref{bound4} with
\begin{displaymath}
  \rho_2 =  \frac{a^2b^4k^2d^3}{q^2}, \quad \alpha_1 = \frac{k}{qr_1r_2^2}, \quad \alpha_2 = \frac{a\sqrt{dn}}{\sqrt{q}}, \quad z = \frac{h}{qr_1r_2^2}, \quad Q = X = q^{\varepsilon}
\end{displaymath}
getting (up to a negligible error)
\begin{displaymath}
\begin{split}
& q^{\varepsilon}  \sum_{abdk} \sum_{n, r_1, r_2}  \frac{w_1(n/N)  w_4(r_1r_2^2/R) r_2}{q^{3/2} (MN)^{1/2}R  a b^2 k d^{3/2}} \\
 &\quad \times \Bigl( \sum_{|h| \leq k} \frac{M^{1/4}q^{3/4}R^{1/2}}{(ak)^{1/2}(dN)^{1/4}}  + \sum_{|h| \geq k}  \frac{qR^{1/2}}{a(dk)^{1/2}N^{1/2}} \Bigl(1+ |h|\min\Bigl(\frac{1}{k},  \frac{\sqrt{M}}{aR\sqrt{dNq} }\Bigr)\Bigr)^{-10}\Bigr).
 \end{split}
\end{displaymath}
The contribution $|h| \leq k$ is 
\begin{displaymath}
 \ll  q^{\varepsilon}  \sum_{abdk} \frac{R^{1/2} N^{1/4}}{q^{3/4} M^{1/4} a^{3/2} k^{1/2} b^2 d^{7/4}} \ll q^{\varepsilon - 9/8} \sum_{k \leq q^{1+\varepsilon}} M^{1/4} \ll q^{\varepsilon - 1/8}
\end{displaymath}
by \eqref{sizes2}. The contribution $|h| \geq k$ is
\begin{displaymath}
\begin{split}
&  \ll q^{\varepsilon}  \sum_{abdk}  \frac{N^{1/2}}{q^{3/2}M^{1/2}  a b^2 k d^{3/2}} \Bigl(  \frac{qR^{1/2}k }{a(dk)^{1/2}N^{1/2}} +  \frac{q^{3/2}R^{3/2}}{M^{1/2}k^{1/2}}\Bigr)\\
 & = q^{\varepsilon} \sum_{abdk} \Bigl(  \frac{R^{1/2}}{q^{1/2} M^{1/2} (abd)^2 k^{1/2}} + \frac{N^{1/2}R^{3/2}}{Mab^2 (kd)^{3/2}}\Bigr) \ll q^{\varepsilon}.
 \end{split} 
\end{displaymath}
by \eqref{sizes2}. This completes the proof.

\end{document}